\documentclass[a4paper]{article}
\usepackage{amssymb,amscd,amsfonts,mathrsfs,amsmath}
\usepackage[all]{xy}
\input amssym.def

\def\double{\mathbb}

\def\cc{{\double C}}
\def\nn{{\double N}}
\def\zz{{\double Z}}

\def\rr{{\double R}}

\newtheorem{theorem}{Theorem}[section]
\newtheorem{lemma}[theorem]{Lemma}

\newtheorem{proposition}[theorem]{Proposition}
\newtheorem{remark}[theorem]{Remark}

\def\Kt{K^{\mathrm{top}}}

\def\cp{\rtimes}

\def\Si{\Sigma}
\def\cinf{C^{\infty}}
\def\cinfc{C^{\infty}_c}

\newcommand{\be}{\begin{equation}}
\newcommand{\ee}{\end{equation}}
\newcommand{\beq}{\begin{eqnarray}}
\newcommand{\eeq}{\end{eqnarray}}
\newcommand{\om}{\omega}
\newcommand{\Om}{\Omega}
\newcommand{\al}{\alpha}
\def\nat{\natural}

\newcommand{\la}{\lambda}
\newcommand{\Ec}{{\mathscr E}}

\newcommand{\Lc}{{\mathscr L}}
\newcommand{\non}{\nonumber}
\newcommand{\eps}{\varepsilon}

\newcommand{\Mc}{{\mathscr M}}

\newcommand{\Ic}{{\mathscr I}}

\def\ch{\mathrm{ch}}

\def\Td{\mathrm{Td}}

\def\tch{\mathrm{\,\slash\!\!\!\! \ch}}

\newcommand{\Tr}{{\mathop{\mathrm{Tr}}}}

\newcommand{\Ac}{{\mathscr A}}

\newcommand{\cqfd}{\hfill\rule{1ex}{1ex}}
\def\deb{\overline{\partial}}

\def\zb{\overline{z}}
\def\ib{\overline{i}}
\def\sb{\overline{s}}
\def\gb{\overline{g}}

\def\wb{\overline{w}}
\def\vb{\overline{v}}

\def\d{\partial}
\def\dd{\mathrm{\bf d}}

\def\hb{\bar{h}}

\def\Bc{{\mathscr B}}

\def\bb{\overline{b}}

\def\hom{{\mathop{\mathrm{Hom}}}}
\def\dom{{\mathop{\mathrm{Dom}}}}
\def\ran{{\mathop{\mathrm{Ran}}}}

\def\hotimes{\hat{\otimes}}

\def\At{\widetilde{A}}

\def\ut{\tilde{u}}
\def\omt{\tilde{\om}}

\def\et{\tilde{e}}

\def\ft{\tilde{f}}
\def\Omh{\widehat{\Omega}}
\def\Th{\widehat{T}}
\def\chih{\widehat{\chi}}
\def\etah{\widehat{\eta}}

\def\Mch{\widehat{\mathscr M}}

\def\Omc{\Omega_c}

\def\mod{\ \mathrm{mod}\ }

\def\supp{\mathrm{supp}\,}

\def\eh{\hat{e}}

\def\uh{\hat{u}}

\def\SL{\mathrm{SL}}

\begin{document}

\begin{center}

{\bf LOCALIZATION OVER COMPLEX-ANALYTIC GROUPOIDS \\
AND CONFORMAL RENORMALIZATION}
\vskip 1cm
{\bf Denis PERROT}
\vskip 0.5cm
Universit\'e de Lyon, Universit\'e Lyon 1,\\
CNRS, UMR 5208 Institut Camille Jordan,\\
43, bd du 11 novembre 1918, 69622 Villeurbanne Cedex, France \\[2mm]
{\tt perrot@math.univ-lyon1.fr}\\[2mm]
\today
\end{center}
\vskip 0.5cm
\begin{abstract}
We present a higher index theorem for a certain class of \'etale one-dimensional complex-analytic groupoids. The novelty is the use of the local anomaly formula established in a previous paper, which represents the bivariant Chern character of a quasihomomorphism as the chiral anomaly associated to a renormalized non-commutative chiral field theory. In the present situation the geometry is non-metric and the corresponding field theory can be renormalized in a purely conformal way, exploiting the complex-analytic structure of the groupoid only. The index formula is automatically localized at the automorphism subset of the groupoid and involves a cap-product with the sum of two different cyclic cocycles over the groupoid algebra. The first cocycle is a trace involving a generalization of the Lefschetz numbers to higher-order fixed points. The second cocycle is a non-commutative Todd class, constructed from the modular automorphism group of the algebra.
\end{abstract}

\vskip 0.5cm


\section{Introduction}

In the previous paper \cite{P6} we presented a general principle allowing to find local representatives for the bivariant Chern character of quasihomomorphisms. The method is based on a \emph{renormalized} bivariant eta-cochain, whose boundary gives the desired local formula. In fact when the quasihomomorphism has even parity this is completely equivalent to the computation of the anomaly associated to a non-commutative chiral gauge theory \cite{P4}, and this explains why the representative of the Chern character is local. The power of this method comes from the considerable freedom in the choice of renormalization: changing the renormalization scheme just amounts to change the local representative of the Chern character. The choice of scheme is however dictated by the geometric situation at hand. In \cite{P6} we illustrated these principles in the ``metric'' situation, i.e. when an abstract Dirac operator is available. In that case the eta-cochain can be renormalized by zeta-function and the local formulas for the bivariant Chern character generalize the Connes-Moscovici formula \cite{CM95}. There exist however non-metric situations, for instance in conformal geometry, where the introduction of a Dirac operator is quite unnatural and more adapted renormalization schemes are needed. The aim of the present paper is to provide a relevant example.\\

We consider a smooth \'etale groupoid $\Gamma$ associated to the action of a discrete group $G$ on the complex plane $\Si$ by local conformal transformations. By this, we mean that each $g\in G$ has a prescribed domain $\dom(g)\subset \Si$ and range $\ran(g)\subset\Si$, and $g$ determines an invertible conformal transformation $\ran(g)\to\dom(g)$. See section \ref{sloc}. A good example is a discrete subgroup of $\SL(2,\cc)$ acting on the plane by homographic transformations. In general the action can be complicated, does not preserve any riemannian metric and may have fixed points. $\Si$ is a Riemann surface, and we let $\deb: \cinfc(\Si)\to \Omc^{0,1}(\Si)$ be the Dolbeault operator acting on compactly supported functions. $\deb$ intertwines the action of $G$ on functions and on differential forms of bidegree $(0,1)$. From the non-commutative geometry viewpoint the groupoid $\Gamma$ is described by the algebraic crossed product $\Ac_0 = \cinfc(\Si)\cp G$. We represent $\Ac_0$ in the algebra of bounded operators on the Hilbert space $H_{\al}$ of square-integrable functions with respect to the measure $(1+|z|)^{\al} d^2z$ over $\Si$, and let $\Bc_0$ be the group ring of $G$. It turns out that for an appropriate choice of weight $\al$, the Dolbeault operator determines a quasihomomorphism of even parity (see \cite{P5})
\be
\rho: \Ac \to \Ec^s \triangleright \Ic^s\hotimes\Bc\ ,
\ee
where $\Ac$ and $\Bc$ are suitable Fr\'echet completions of $\Ac_0$ and $\Bc_0$, the algebra $\Ic$ is the Schatten ideal of $p$-summable operators $\Lc^p(H_{\al})$ for any choice of $p>2$, $\Ec=\Lc(H_{\al})\hotimes\Bc$, and $\hotimes$ is the projective tensor product. The Riemann-Roch-Grothendieck theorem stated in \cite{P5} establishes the compatibility between the pushforward maps induced by $\rho$ on various versions of $K$-theory and cyclic homology, including secondary characteristic classes. We will focus on topological $K$-theory and periodic cyclic homology only, i.e. the homotopy invariants of $\Ac$ and $\Bc$. In that case the theorem reduces to a commutative diagram:
\be
\vcenter{\xymatrix{
\Kt_i(\Ic\hotimes\Ac) \ar[r]^{\rho_!} \ar[d] \ar@{.>}[rd]  & \Kt_i(\Ic\hotimes\Bc) \ar[d] \\
HP_i(\Ac) \ar[r]_{\ch(\rho)} & HP_i(\Bc) }}  \qquad i\in \zz_2\ . \label{dg}
\ee
The vertical arrow $\Kt_i(\Ic\hotimes\cdot)\to HP_i(\cdot)$ is the Chern character in periodic cyclic homology and $\rho_!$ is the pushforward map in topological $K$-theory. The bottom arrow is the bivariant periodic Chern character $\ch(\rho)\in HP^0(\Ac,\Bc)$. We are interested in the explicit calculation of the diagonal map (dashed arrow): from cyclic homology invariants of $\Bc$ one thus gets $K$-theoretic invariants of $\Ac$. This is useful for example in the formulation of higher index theorems (see e.g. \cite{CM90, P3}). Note that the methods developed for the equivariant index theorem of \cite{P3} cannot be applied here because the action of $G$ does not preserve any Riemannian metric on $\Si$. For this reason heat kernel or zeta-function renormalizations are not suited and in order to obtain local formulas we shall exploit the \emph{conformal geometry} only, i.e. the complex-analytic structure of $\Gamma$. Alternatively, in the odd case $i=1$ the diagonal map is exactly computed by a chiral anomaly \cite{P6} and the index theorem is thus reduced to the renormalization of a chiral conformal field theory over $\Gamma$. The even case $i=0$ is then covered by Bott periodicity. \\

As stressed in \cite{P6} renormalization requires to work with the dense subalgebra $\Ac_0\subset\Ac$. Indeed the cyclic cocycles obtained are distributions and make sense on the algebra of smooth functions with compact support $\cinfc(\Si)\cp G$. We shall not address here the question of extending these cocycles to the Fr\'echet algebra $\Ac$. Hence let us choose an even $K$-theory class $[e] \in \Kt_0(\Ic\hotimes\Ac)$ represented by an idempotent $e\in M_{\infty}(\Ac_0)$, and an odd class $[u] \in \Kt_1(\Ic\hotimes\Ac)$ represented by an invertible $u \in M_{\infty}(\Ac_0)^+$ such that $u-1 \in M_{\infty}(\Ac_0)$ (the symbol $^+$ denotes unitalization). Their respective images under the canonical homomorphism $\tilde{\rho}: \Ac_0 \to \Ac_0\otimes\Bc$ have Chern characters 
\be
\ch(\tilde{\rho}(e)) \in HP_0(\Ac_0\otimes\Bc)\ ,\qquad \ch(\tilde{\rho}(u)) \in HP_1(\Ac_0\otimes\Bc)\ .
\ee
(Here the algebraic tensor product $\Ac_0\otimes \Bc$ is treated as a discrete algebra). The diagonal map of (\ref{dg}) then corresponds to the cap-product of these Chern characters with an appropriate cyclic cocycle $\varphi$ of even degree over $\Ac_0$:
\be
\varphi \cap : HP_i(\Ac_0\otimes\Bc) \to HP_i(\Bc)\ .
\ee
The cocycle $\varphi$ emerges from the anomaly formula. For this reason it is automatically localized at the automorphism subset of the groupoid $\Gamma$, i.e. the pairs $(g,z_0)\in G\times \Si$ given by a fixed point $z_0$ for a conformal mapping $g$. The set of automorphisms is the union of the discrete set $\Gamma_f$ of isolated automorphisms, and the 1-dimensional complex manifold $\Gamma_{\infty}$ of non-isolated automorphisms (which contains the set of units $\Si$). An isolated automorphism $(g,z_0)$ has an order $n<\infty$ corresponding to the first non-vanishing term of the Taylor expansion of the holomorphic function $g(z)-z$ around $z_0$. The non-isolated automorphisms have order $n=\infty$ by definition ($g(z)=z$ around $z_0$). Hence $\varphi$ splits as the sum of two cocycles. The first part is a cyclic cocycle of degree zero, i.e. a trace $\Phi(\Gamma):\Ac_0 \to \cc$ localized at isolated automorphisms:
\be
\Phi(\Gamma)(a) = \sum_{ (g,z_0) \in\Gamma_f } \frac{-1}{(n-1)!}\, \d_z^{n-1} \big( H^n_{g,z_0}(z)\,a(g,z) \big)_{z=z_0}\qquad \forall a\in \Ac_0 \label{fix}
\ee
where $H^n_{g,z_0}(z)$ is the holomorphic function $(z-z_0)^n/(g(z)-z)$ and $n$ is the order of $(g,z_0)$. The second part is a cyclic cocycle of degree two localized at the manifold of non-isolated automorphisms. Consider the crossed product algebra of differential forms with compact support $\Omc^*(\Si)\cp G$. It contains $\Ac_0$ as subalgebra of degree zero and is gifted with the de Rham differential $d=\d + \deb$ acting on $\Omc^*(\Si)$ only and commuting with $G$. There is also a non-commutative differential $\delta$ coming from the \emph{modular derivative} $D$ over $\Ac_0$:
\be
\delta = [\d, D]\ ,\qquad \delta^2=0\ .
\ee
$D$ is the generator of the modular automorphism group of $\Ac_0$ associated to the smooth euclidian volume form $d\zb\wedge dz/2i$ over $\Si$ (see \cite{P0}). From these ingredients we form the differential $\nabla=d - \frac{1}{2}\delta$, $\nabla^2=0$ and define the Todd class of $\Gamma$ as the cyclic 2-cocycle
\be
\Td(\Gamma)(a_0,a_1,a_2) = \int_{\Gamma_{\infty}}a_0\nabla a_1\nabla a_2\ ,\qquad \forall a_i\in \Ac_0\ .
\ee
It was shown in \cite{P0} how the classical counterpart of the modular differential $\delta$ accounts for the curvature of a Riemann surface endowed with a K\"ahler metric, and thus explains why $\Td(\Gamma)$ is the correct generalization of the usual Todd class. Note by the way that $\delta$ was introduced in \cite{P0} as one of the generators of the Connes-Moscovici Hopf algebra \cite{CM98}. Using the anomaly formula one gets for free:

\begin{theorem}
Let $[e]\in \Kt_0(\Ic\hotimes\Ac)$ and $[u]\in \Kt_1(\Ic\hotimes\Ac)$ be topological $K$-theory classes represented by an idempotent $e\in M_{\infty}(\Ac_0)$ and an invertible element $u\in M_{\infty}(\Ac_0)^+$ respectively. Their images under the diagonal map of the commutative diagram (\ref{dg}) are obtained by cap product of the Chern characters $\ch(\tilde{\rho}(e))\in HP_0(\Ac_0\otimes\Bc)$ and $\ch(\tilde{\rho}(u))\in HP_1(\Ac_0\otimes\Bc)$ with the cyclic cocycle over $\Ac_0$:
\be
\varphi = \Phi(\Gamma) + \Td(\Gamma)\ .
\ee
\end{theorem}

See Theorem \ref{tcup} for an explicit formula. The cap-product with the Todd class is not very surprising, since it is the general form expected in a non-commutative Riemann-Roch-Grothendieck theorem (see the work of Brodzki, Mathai, Rosenberg and Szabo \cite{BMRS}), but the trace $\Phi$ is more exotic. It is not obvious from (\ref{fix}) but nevertheless true that $\Phi$ is invariantly defined, i.e. only depends on the complex-analytic structure of $\Gamma$ and not on the particular choice of complex coordinate system $z$ (Lemma \ref{linv}). The contribution of an isolated automorphism $(g,z_0)$ of order $n$ can be computed in low degrees,
\beq
\lefteqn{\frac{-1}{(n-1)!}\, \d_z^{n-1}\big(H_{g,z_0}^n(z) a(g,z)\big)_{z=z_0} = } \non\\
&(n=1)&:\qquad \frac{1}{1-g'(z_0)}\, a(z_0)    \non\\
&(n=2)&:\qquad \frac{2}{g''(z_0)} \left( \frac{1}{3} \frac{g'''(z_0)}{g''(z_0)} \, a(z_0) -  \d_z a(z_0) \right)   \non\\
&(n=3)&:\qquad \frac{3}{2g'''(z_0)} \left( \frac{1}{10} \frac{g^{(5)}(z_0)}{g'''(z_0)}\right.  a(z_0) - \frac{1}{8} \left(\frac{g^{(4)}(z_0)}{g'''(z_0)}\right)^2 \, a(z_0)  \non\\
&& \qquad \qquad \qquad \qquad \left. + \frac{1}{2} \frac{g^{(4)}(z_0)}{g'''(z_0)} \, \d_z a(z_0) -  \d_z^2 a(z_0) \right)   \non
\eeq
and depends on the jets of the conformal mapping $g$ up to order $2n-1$. For $n=1$ one recognizes the Lefschetz number $1/(1-g'(z_0))$ and it pairs non-trivially with the $K$-theory of $\Ac_0$. For $n\geq 2$ one gets new traces over $\Ac_0$. They do not appear in usual fixed-point theorems because the latter deal with \emph{isometries}, which have only fixed points of order one. I don't know if these traces pair non-trivially with $K$-theory, or even if their periodic cyclic cohomology class is non-zero. One can expect that their evaluation on $K$-theory classes in the range of an ``assembly map'' (whatever it is) should vanish because a conformal transformation $g$ associated to a fixed point $z_0$ of order $\geq 2$ cannot generate a finite group. Nevertheless, these evanescent terms are traces over $\Ac_0$ in their own right and may be of independent interest.\\

The paper is organized as follows. In section \ref{sloc} we construct the quasihomomorphism associated to the action of $G$ on $\Si$ and recall from \cite{P6} how local representatives of its bivariant Chern character are related to the chiral anomaly of a non-commutative conformal field theory. Then in section \ref{sren} we explicitly renormalize the conformal field theory and compute the anomaly, spontaneously localized at the fixed points of conformal mappings. These calculations are rather straightforward but involve distributions and have to be performed with care. Finally the index theorem is stated in section \ref{sind}.

\section{Local anomaly formula}\label{sloc}

Let $\Si=\cc$ be the complex plane. We regard $\Si$ as a Riemann surface with its canonical orientation given by the complex structure. In the following $z$ denotes the canonical complex coordinate system on $\Si$. Abusively $z$ will sometimes also denote a point in $\Si$. Let $\Omc^*(\Si)$ be the differential graded algebra of smooth, complex-valued  differential forms with compact support. The degree zero subalgebra $\Omc^0(\Si)$ coincides with the algebra of smooth compactly supported functions $\cinfc(\Si)$. Any one-form $A\in\Omc^1(\Si)$ can be decomposed as a sum $A=dz A_z+d\zb A_{\zb}$, and we let $\Omc^{1,0}(\Si)$ and $\Omc^{0,1}(\Si)$ be the subspaces of one-forms proportional to $dz$ and $d\zb$ respectively. The partial derivatives with respect to $z$ and $\zb$ will be denoted by $\d_z$ and $\d_{\zb}$, hence the de Rham differential reads $d=dz\d_z+d\zb\d_{\zb}$. The \emph{Dolbeault operator} is defined as usual as the part proportional to $d\zb$:
\be
\deb=d\zb\d_{\zb} : \cinfc(\Si)\to \Omc^{0,1}(\Si)\ .
\ee
For notational consistency with \cite{P6} the symbol $Q$ will be used instead of $\deb$. We denote improperly by $Q^{-1}$ the associated Green's operator. Thus $Q^{-1}$ is a linear map from $\Omc^{0,1}(\Si)$ to the (non-compactly supported) smooth functions $\cinf(\Si)$, given by the integral
\be
(Q^{-1}\cdot A)(w)=\int_{\Si} d^2z \, \frac{A_{\zb}(z)}{\pi(w-z)}\ ,
\ee
for all $A=d\zb A_{\zb} \in \Omc^{0,1}(\Si)$ at any point $w\in \Si$. Here $d^2z=d\zb\wedge dz/2i$ is the euclidean volume form, and the Cauchy kernel $1/(w-z)$ is a distribution over $\Si\times \Si$ with singular support on the diagonal. The notation $Q^{-1}$ is justified by the fact that for any $f\in \cinfc(\Si)$ and $A\in \Omc^{0,1}(\Si)$ one has $Q^{-1}Q\cdot f=f$ and $QQ^{-1}\cdot A=A$. Indeed this can be easily verified using the distributional relation
$$
\d_{\zb}\, \frac{1}{z-w}=\pi \delta^2(z-w)\ ,
$$
where $\delta^2(z)$ is the Dirac mass at $z=0$. We would like to complete $\cinfc(\Si)$ into a Hilbert space of square-integrable functions with respect to some measure: for any weight $\al\in\rr$ endow $\cinfc(\Si)$ with the norm
\be
\|\xi\|_{\al}=\Big(\int_{\Si} d^2z\, (1+|z|)^{\al} |\xi(z)|^2 \Big)^{1/2}\ ,\quad \forall \xi\in\cinfc(\Si)\ .
\ee
Let $H_{\al}$ be the Hilbert space completion of $\cinfc(\Si)$ with respect to this norm, and observe that any one-form $A\in \Omc^{0,1}(\Si)$ acts by pointwise multiplication on the space of smooth functions and thus defines a linear map $\cinfc(\Si)\to \Omc^{0,1}(\Si)$. 

\begin{lemma}\label{lcom}
For any one-form $A\in\Omc^{0,1}(\Si)$, the composite operator $Q^{-1}A:\cinfc(\Si)\to\cinf(\Si)$ extends to a compact operator on $H_{\al}$ whenever $\al<-1$, and more precisely $Q^{-1}A$ is in the Schatten class $\Lc^p(H_{\al})$ for any $p>2$.
\end{lemma}
{\it Proof:} We will show that the formal adjoint of $Q^{-1}A$ extends to a compact operator. The Hilbert norm $\| \cdot \|_{\al}$ comes from the inner product
$$
\langle \xi, \eta \rangle_{\al} = \int d^2z\, (1+|z|)^{\al} \overline{\xi(z)}\, \eta(z)
$$
for all $\xi,\eta \in \cinfc(\Si)$. The adjoint $(Q^{-1}A)^*$ with respect to this product reads 
$$
(Q^{-1}A)^* \xi(z) = (1+|z|)^{-\al}\overline{A_{\zb}(z)} \int d^2w\, \frac{(1+|w|)^{\al}}{\pi(\wb-\zb)}\xi(w)\ .
$$
Choose $\al < -1$ and let us show that $(Q^{-1}A)^*$ extends to a bounded map from $H_{\al}$ to the Sobolev space $W_1$ gifted with the norm
$$
\|\xi\|_W = \Big( \int d^2z\, (|\xi(z)|^2 + |\d_z\xi(z)|^2) \Big)^{1/2}\ .
$$
Define the function with compact support $f(z)=(1+|z|)^{-\al}\overline{A_{\zb}(z)}$. Using the identity $\d_z \frac{1}{\zb-\wb} = \pi \delta^2(z-w)$ one gets
\beq
\lefteqn{ (\| (Q^{-1}A)^* \xi(z) \|_W)^2 \leq \int d^2z\, \left|f(z)\int d^2w\, \frac{(1+|w|)^{\al}}{\pi(\wb-\zb)}\xi(w) \right|^2 + } \non\\
&& \mspace{100mu} \int d^2z\, \left|\d_z f(z)\int d^2w\, \frac{(1+|w|)^{\al}}{\pi(\wb-\zb)}\xi(w) - f(z)(1+|z|)^{\al}\xi(z) \right|^2 \non
\eeq
We estimate the first term of the r.h.s. as follows:
\beq
\lefteqn{\int d^2z\, \left| f(z)\int d^2w\, \frac{(1+|w|)^{\al}}{\pi(\wb-\zb)}\xi(w) \right|^2 } \non\\
&&\leq \int d^2z d^2w_1 d^2w_2\, |f(z)|\frac{(1+|w_1|)^{\al}}{\pi|w_1-z|} |\xi(w_1)| \, |f(z)|\frac{(1+|w_2|)^{\al}}{\pi|w_2-z|} |\xi(w_2)|  \non \\
&& \leq \int d^2z d^2w_1 d^2w_2\, |f(z)|^2\frac{(1+|w_1|)^{\al}}{\pi|w_1-z|} \, \frac{(1+|w_2|)^{\al}}{\pi|w_2-z|} |\xi(w_1)|^2 \ , \non
\eeq
where the last step uses the Cauchy-Schwarz inequality. Since $\al < -1$ the integral $I(z)=\int d^2w_2\,\frac{(1+|w_2|)^{\al}}{|w_2-z|} $ converges to a continuous function of $z$. Afterwards, the integral $\int d^2z |f(z)|^2\frac{I(z)}{|w_1-z|} $ is uniformly bounded with respect to $w_1$. Hence there is a constant $C_{\al}(f)$ such that 
\beq
\int d^2z\, \left|f(z)\int d^2w\, \frac{(1+|w|)^{\al}}{\pi(\wb-\zb)}\xi(w) \right|^2 &\leq& C_{\al}(f)\, \int d^2w_1 (1+|w_1|)^{\al}|\xi(w_1)|^2 \non\\
&\leq& C_{\al}(f)\, (\|\xi\|_{\al})^2 \non
\eeq
To estimate the second term we develop 
\beq
\lefteqn{\int d^2z \left|\d_z f(z)\int d^2w\, \frac{(1+|w|)^{\al}}{\pi(\wb-\zb)}\xi(w) - f(z)(1+|z|)^{\al}\xi(z) \right|^2 } \non \\
&\leq& \int d^2z \left|\d_z f(z)\int d^2w\, \frac{(1+|w|)^{\al}}{\pi(\wb-\zb)}\xi(w)\right|^2 + \int d^2z \left|f(z)(1+|z|)^{\al}\xi(z) \right|^2 \non\\
&& + 2\int d^2z d^2w\, |\d_z f(z)| \frac{(1+|w|)^{\al}}{\pi|w-z|}|\xi(w)| \, |f(z)|(1+|z|)^{\al}|\xi(z)| \non \\
&\leq& C_{\al}(\d_zf)(\|\xi\|_{\al})^2  + C' \int d^2z\, (1+|z|)^{\al} |\xi(z)|^2 \non\\
&& + 2 \int d^2z d^2w\, |\tilde{f}(z)|(1+|z|)^{\al}  \frac{(1+|w|)^{\al}}{\pi|w-z|}|\xi(w)||\xi(z)|\ , \non
\eeq
where $\tilde{f}(z)$ is the compactly supported function $\d_z f(z) f(z)$ and $C'$ is the supremum of $|f(z)|(1+|z|)^{\al}$. The Cauchy-Schwarz inequality implies
\beq
\lefteqn{\int d^2z d^2w\, |\tilde{f}(z)|(1+|z|)^{\al}  \frac{(1+|w|)^{\al}}{|w-z|}|\xi(w)||\xi(z)| } \non\\
&&\qquad \leq  \Big( \int d^2z d^2w\, |\tilde{f}(z)| \frac{(1+|z|)^{\al}}{|w-z|} (1+|w|)^{\al} |\xi(w)|^2 \Big)^{1/2} \non\\
&& \qquad\qquad  \times \Big( \int d^2z d^2w\, |\tilde{f}(z)| \frac{(1+|w|)^{\al}}{|w-z|}(1+|z|)^{\al} |\xi(z)|^2 \Big)^{1/2}\ . \non
\eeq   
The integral $\int d^2z |\tilde{f}(z)| \frac{(1+|z|)^{\al}}{|w-z|} $ is uniformly bounded with respect to $w$, and $\int d^2w |\tilde{f}(z)| \frac{(1+|w|)^{\al}}{|w-z|}$ is uniformly bounded with respect to $z$. We conclude that 
$$
\| (Q^{-1}A)^* \xi \|_W \leq C''_{\al}(f) \, \|\xi\|_{\al}
$$
for some constant $C''_{\al}(f)$, and the adjoint $(Q^{-1}A)^*$ extends to a bounded map from $H_{\al}$ to the Sobolev space $W_1$. Moreover its range is a set of functions with support contained in the compact support of the one-form $A$. It follows from Rellich's lemma that $(Q^{-1}A)^*$ is a compact operator on $H_{\al}$, and more precisely is $p$-summable for any $p>2$.  \cqfd \\

Let $G$ be a discrete group acting by local conformal transformations on $\Si$. It means that any $g\in G$ has prescribed domain $\dom(g)$ and range $\ran(g)$ which are (possibly empty) open subsets of $\Si$, and $g$ determines an invertible conformal mapping $\dom(g)\to\ran(g)$ subject to the following conditions:
\begin{itemize}
\item The unit $1\in G$ acts by the identity on $\dom(1)\subset\Si$.
\item $\dom(g^{-1})=\ran(g)$ for any $g\in G$ and $g^{-1}$ determines the inverse mapping of $g$.
\item $\dom(gh) \supset h^{-1}(\dom(g))\cap\dom(h)$ for any $g,h\in G$. 
\end{itemize}
In particular one necessarily has $\dom(g)\subset\dom(1)$. We don't need to impose the equality $\dom(1)=\Si$. Remark that in general the conformal mapping $\dom(g)\to\ran(g)$ induced by $g$ does \emph{not} specify $g$ as an element of the group $G$: for example, $G$ is an arbitrary group acting by the identity on $\Si$ with $\dom(g)=\Si$ for any $g\in G$. Nevertheless, for convenience we will often refer to ``the conformal mapping $g:\dom(g)\to\ran(g)$'' instead of ``the conformal mapping $\dom(g)\to\ran(g)$ induced by $g\in G$''. This is abusive but should not create too much confusion. To this data we associate the crossed product 
\be
\Ac_0=\cinfc(\Si)\cp G\ .
\ee
It is the algebra linearly generated by finite sums of symbols $fU^*_g$, with $g\in G$ and $f\in \cinfc(\dom(g))$. The product is given by convolution:
$$
(f_1U^*_{g_1})(f_2U^*_{g_2}) = f_1 f_2^{g_1}U^*_{g_2g_1}\ ,
$$
where $f_2^{g_1}$ denotes the pullback of the function $f_2$ by the mapping induced by $g_1$. The conditions on the domains imply that the product is well-defined and associative. We also consider the group ring of $G$, as the algebra $\Bc_0$ linearly generated by finite sums of symbols $U^*_g$, with product $U^*_{g_1}U^*_{g_2} = U^*_{g_2g_1} $. The point is that the Dolbeault operator $Q$ gives rise to a quasihomomorphism between adequate Fr\'echet completions of the algebras $\Ac_0$ and $\Bc_0$. Firstly, $\Ac_0$ is canonically represented in the algebra of bounded operators $\Lc(H_{\al})$. Indeed, any element $fU^*_g\in \Ac_0$ acts as a linear operator $f r(g)_+:\cinfc(\Si)\to\cinfc(\Si)$, with $r(g)_+$ the action of the conformal mapping $g$ on smooth functions by pullback:
\be
(fr(g)_+\cdot \xi)(z)=f(z)\xi(g(z))\ ,\qquad \forall \ \xi\in \cinfc(\Si)\ ,
\ee
and this makes sense because $\supp(f)\subset \dom(g)$. Moreover, since $f$ has compact support $fr(g)_+$ extends to a bounded operator on $H_{\al}$, and one easily checks that it yields an algebra homomorphism $\Ac_0\to\Lc(H_{\al})$. In the same manner, $fU^*_g$ defines a linear operator $fr(g)_-:\Omc^{0,1}(\Si)\to \Omc^{0,1}(\Si)$, where $r(g)_-$ is the action of the conformal mapping $g$ on one-forms by pullback. It turns out that conjugation with $Q$ also gives a bounded operator $Q^{-1}f r(g)_- Q$ on $H_{\al}$ provided $\al<-1$. It suffices to remark that the Dolbeault operator commutes with conformal transformations, hence 
$$
Q^{-1}fr(g)_-Q = Q^{-1}fQr(g)_+ = fr(g)_+ - Q^{-1}[Q,f]r(g)_+ \ .
$$
The commutator $[Q,f]$ coincides with multiplication by the one-form $Q\cdot f=d\zb\d_{\zb}f \in \Omc^{0,1}(\Si)$, therefore $Q^{-1}[Q,f]r(g)_+$ is compact by Lemma \ref{lcom}. Hence $Q^{-1}fr(g)_-Q$ yields another representation $\Ac_0\to \Lc(H_{\al})$, which differs from the former by compact operators. Consequently, we obtain two algebra homomorphisms $(\rho_+,\rho_-): \Ac_0\to \Lc(H_{\al})\otimes \Bc_0$ by setting
\be
\rho(fU^*_g)_+= fr(g)_+\otimes U^*_g\ ,\qquad \rho(fU^*_g)_-= Q^{-1}fr(g)_-Q \otimes U^*_g\ ,
\ee
and the difference $\rho_+-\rho_-$ is a linear map from $\Ac_0$ to the two-sided ideal $\Lc^p(H_{\al})\otimes\Bc_0$ for any $p>2$. To get a quasihomomorphism we need to complete: choose a Fr\'echet $m$-algebra $\Bc$ which contains $\Bc_0$ as a dense subalgebra. By viewing $\Lc(H_{\al})$ as a Banach algebra with operator norm, the completed projective tensor product $\Ec=\Lc(H_{\al})\hotimes \Bc$ is a Fr\'echet $m$-algebra completion of $\Lc(H_{\al})\otimes \Bc_0$. Let $\Ic$ be the Banach algebra $\Lc^p(H_{\al})$ endowed with the Schatten norm $\|\cdot\|_p$. The completed projective tensor product $\Ic\hotimes\Bc$ is a Fr\'echet $m$-algebra completion of $\Lc^p(H_{\al})\otimes\Bc_0$. The continuous inclusion $\Ic\to\Lc(H_{\al})$ induces a continuous homomorphism $\Ic\hotimes\Bc\to\Ec$. In the sequel \emph{we assume that the latter is injective}, which promotes $\Ic\hotimes\Bc$ to a (non-closed) two-sided ideal in $\Ec$. This assumption is not really crucial but it simplifies the discussion. Following \cite{P5} we introduce the semi-direct sum 
\be
\Ec^s_+=\Ec\ltimes \Ic\hotimes\Bc\ .
\ee
It is $\Ec\oplus \Ic\hotimes \Bc$ as a locally convex vector space, and the product puts as many terms as possible in $\Ic\hotimes\Bc$. Then $\Ec^s_+$ is a Fr\'echet $m$-algebra containing $\Ic\hotimes\Bc$ as a closed two-sided ideal. This situation is depicted by $\Ec\triangleright \Ic\hotimes\Bc$. The algebra $\Ec^s_+$ is endowed with a canonical action of the group $\zz_2$ by automorphisms, and we define $\Ec^s$ as the $\zz_2$-graded crossed product algebra $\Ec^s_+\cp\zz_2$. It contains $\Ic^s\hotimes\Bc$ as a (not necessarily closed) two-sided ideal where $\Ic^s$ is the $\zz_2$-graded algebra of $2\times 2$ matrices over $\Ic$, and one has $\Ec^s\triangleright\Ic^s\hotimes\Bc$. We refer to \cite{P5} for details.  \\
Now the pair $(\rho_+,\rho_-):\Ac_0\to \Ec$ induces a homomorphism $\rho:\Ac_0\to \Ec^s_+$ by setting $\rho=\rho_-\oplus(\rho_+-\rho_-)$. In fact $\rho$ is injective, hence we may complete $\Ac_0$ into a Fr\'echet $m$-algebra by taking its closure in $\Ec^s_+$. We should keep in mind that $\Ac$ depends on the choice of completion $\Bc$ and of weight $\al$ for the Hilbert space $H_{\al}$, but the renormalization calculations performed afterwards will entirely deal with $\Ac_0$ and hence are completely independent of these choices. To summarize one gets a $p$-summable quasihomomorphism of even parity
\be
\rho: \Ac \to \Ec^s\triangleright \Ic^s\hotimes\Bc
\ee
for any $p>2$, according to Definition 3.1 of \cite{P5}. $\rho$ has a Chern character $\ch(\rho)\in HP^0(\Ac,\Bc)$ in bivariant periodic cyclic cohomology of even degree. Its compatibility with the pushforward map $\rho_!$ in topological  $K$-theory is expressed by the Riemann-Roch-Grothendieck theorem of \cite{P5}, through the commutative diagram
\be
\vcenter{\xymatrix{
\Kt_i(\Ic\hotimes\Ac) \ar[r]^{\rho_!} \ar[d] & \Kt_i(\Ic\hotimes\Bc) \ar[d] \\
HP_i(\Ac) \ar[r]^{\ch(\rho)} & HP_i(\Bc) }}  \label{rrg}
\ee
where $\Kt_i(\Ic\hotimes\cdot)\to HP_i(\cdot)$ is the Chern character in periodic cyclic homology and $i\in\zz_2$. Since the quasihomomorphism is $p$-summable for any $p>2$, $\ch(\rho)$ is represented by a hierarchy of non-periodic bivariant Chern characters $\ch^n(\rho)\in HC^n(\Ac,\Bc)$ in all even degrees $n\geq 2$. They are related by the $S$-operation in bivariant cyclic cohomology: $\ch^{n+2}(\rho)\equiv S \ch^n(\rho) \in HC^{n+2}(\Ac,\Bc)$. Hence all the $\ch^n(\rho)$'s induce the same map $HP_i(\Ac)\to HP_i(\Bc)$. The construction of bivariant Chern characters requires the choice of quasi-free extensions; in the present situation we work with the universal free extensions 
$$
0 \to J\Ac \to T\Ac \to \Ac \to 0\ ,\qquad 0 \to J\Bc \to T\Bc \to \Bc \to 0 \ ,
$$
where $T\Ac$ is the tensor $m$-algebra over $\Ac$ and its ideal $J\Ac$ is the kernel of the multiplication map $T\Ac\to\Ac$. According to the Cuntz-Quillen formalism \cite{CQ1}, the pro-algebra $\Th\Ac=\varprojlim_k T\Ac/(J\Ac)^k$ calculates the periodic cyclic homology of $\Ac$ as the homology of the $X$-complex
\be
X(\Th\Ac):\quad \Th\Ac\ \mathop{\rightleftarrows}^{\nat \dd}_{\bb}\  \Om^1\Th\Ac_{\nat}\ ,
\ee
where $\Om^1\Th\Ac_{\nat}$ is the universal $\Th\Ac$-bimodule of noncommutative one-forms over $\Th\Ac$ quotiented by the commutator subspace $[\Th\Ac,\Om^1\Th\Ac]$. The map $\nat\dd$ is the universal derivation $\dd$ followed by the quotient map $\nat$, and $\bb$ is the Hochschild boundary $\nat x\dd y\to [x,y]$. The bivariant Chern character $\ch^n(\rho)$ is realized as a bivariant cocycle in the Hom-complex $\hom(X(\Th\Ac),X(\Th\Bc))$, i.e. a chain map. We first lift the quasihomomorphism $\rho:\Ac\to \Ec^s\triangleright \Ic^s\hotimes\Bc$ to a quasihomomorphism
$$
\rho_*: T\Ac \to \Mc^s\triangleright \Ic^s\hotimes T\Bc\ ,
$$
where $\Mc=\Lc(H_{\al})\hotimes T\Bc$ (again we have to suppose that $\Ic\hotimes T\Bc \to \Mc$ is injective). The homomorphism $\rho_*$ is constructed from a pair of homomorphisms $(\rho_{*+},\rho_{*-}):T\Ac\to \Mc$ defined on a chain $x= f_1U^*_{g_1} \otimes\ldots\otimes f_kU^*_{g_k} \in T\Ac$ by the formulas
\beq
\rho_{*}(x)_+ &=& \big(f_1f_2^{g_1} \ldots f_k^{g_{k-1}\ldots g_1}r(g_k\ldots g_1)_+\big) \otimes (U^*_{g_1}\otimes \ldots \otimes U^*_{g_k})\ , \non\\
\rho_{*}(x)_- &=& Q^{-1}\big(f_1f_2^{g_1} \ldots f_k^{g_{k-1}\ldots g_1}r(g_k\ldots g_1)_-\big)Q \otimes (U^*_{g_1}\otimes \ldots \otimes U^*_{g_k})\ . \non
\eeq
We use the convenient representation of the $\zz_2$-graded algebra $\Mc^s$ in terms of $2\times 2$ matrices \cite{P5}, and introduce the odd multiplier $F=\bigl( \begin{smallmatrix}
0 & 1 \\
1 & 0 \end{smallmatrix} \bigr)$. Then we can write $\rho_*(x)= \bigl( \begin{smallmatrix}
\rho_*(x)_+ & 0 \\
0 & \rho_*(x)_- \end{smallmatrix} \bigr)\in \Mc^s_+$ for all $x\in T\Ac$. Since the difference $\rho_*(x)_+-\rho_*(x)_-$ is compact one has $[F,\Mc^s_+]\subset\Ic^s\hotimes T\Bc$. The homomorphism $\rho_*$ is compatible with the adic filtrations induced by the ideals $J\Ac$, $J\Bc$ and thus extends to a homomorphism of pro-algebras $\rho_*:\Th\Ac\to \Mch^s_+$. The bivariant Chern character of degree $n$ is given by the composition of chain maps
\be
\ch^n(\rho) : X(\Th\Ac) \stackrel{\gamma}{\longrightarrow} \Omh\Th\Ac \stackrel{\rho_*}{\longrightarrow} \Omh\Mch^s_+ \stackrel{\chih^n}{\longrightarrow} X(\Th\Bc) \ ,
\ee
where $\gamma$ is the homotopy equivalence between the $X$-complex and the completion of Connes' $(b+B)$-complex \cite{C1} of noncommutative differential forms $\Omh$ (see for example \cite{P5}), and $\chih^n$ has only two non-zero components $\chih^n_0:\Om^n \Mch^s_+\to \Th\Bc$ and $\chih^n_1:\Om^{n+1}\Mch^s_+\to \Om^1\Th\Bc_{\nat}$ given by
$$
\chih^n_0(x_0\dd x_1\ldots\dd x_n) = (-)^n\frac{\Gamma(1+\frac{n}{2})}{(n+1)!} \sum_{\la\in S_{n+1}} \eps(\la)\, \tau(x_{\la(0)}[F,x_{\la(1)}]\ldots [F,x_{\la(n)}])
$$
$$
\chih^n_1(x_0\dd x_1\ldots\dd x_{n+1}) = (-)^n\frac{\Gamma(1+\frac{n}{2})}{(n+1)!} \sum_{i=1}^{n+1}  \tau\nat(x_0[F,x_1]\ldots\dd x_i \ldots [F,x_{n+1}])\ .
$$
Here $S_{n+1}$ is the cyclic permutation group of $n+1$ elements, $\eps$ is the signature of permutations (here it is always 1 since $n$ is even), and $\tau$ is the operator supertrace. For $n=2$ it is necessary to replace $\tau$ by $\tau'=\frac{1}{2}\tau(F[F,\ ])$ in order to ensure traceability. The relation $\ch^{n+2}(\rho)\equiv S\ch^n(\rho)$ is a consequence of the fact that the difference between the cocycles $\chih^n$ and $\chih^{n+2}$ is a coboundary in the $\zz_2$-graded complex $\hom(\Omh\Mch^s_+,X(\Th\Bc))$:
$$
\chih^n-\chih^{n+2}=(\nat\dd\oplus\bb)\circ \etah^{n+1} + \etah^{n+1}\circ (b+B)\ .
$$
The eta-cochain $\etah^{n+1}$ is an odd element of the above Hom-complex and has only two non-zero components $\etah^{n+1}_0:\Om^{n+1}\Mch^s_+\to \Th\Bc$ and $\etah^{n+1}_1:\Om^{n+2}\Mch^s_+\to \Om^1\Th\Bc_{\nat}$ given by
\beq
\lefteqn{\etah^{n+1}_0(x_0\dd x_1\ldots\dd x_{n+1}) =  \frac{\Gamma(\frac{n}{2}+1)}{(n+2)!} \, \frac{1}{2}\tau\Big (F x_0[F,x_1]\ldots [F,x_{n+1}]+ } \non \\
&&\qquad \qquad \qquad  \sum_{i=1}^{n+1}(-)^{(n+1)i} [F,x_i]\ldots [F,x_{n+1}] Fx_0 [F,x_1]\ldots [F,x_{i-1}] \Big) \non
\eeq
\beq
\lefteqn{\etah^{n+1}_1(x_0\dd x_1\ldots\dd x_{n+2}) =} \non\\
&& \frac{\Gamma(\frac{n}{2}+1)}{(n+3)!} \sum_{i=1}^{n+2}  \frac{1}{2}\tau\nat(ix_0 F + (n+3-i)Fx_0)[F,x_1]\ldots\dd x_i \ldots [F,x_{n+2}]\ .\non
\eeq
Define $\tch^{n+1}(\rho)$ as the composition $X(\Th\Ac) \stackrel{\gamma}{\to} \Omh\Th\Ac \stackrel{\rho_*}{\to} \Omh\Mch^s_+ \stackrel{\etah^{n+1}}{\to} X(\Th\Bc)$. Hence the transgression relation $\ch^n(\rho)-\ch^{n+2}(\rho)=[\d,\tch^{n+1}(\rho)]$ holds in the complex $\hom(X(\Th\Ac),X(\Th\Bc))$, where $\d$ denotes the $X$-complex boundary and $[\ ,\ ]$ is the graded commutator.\\

These formulas are not really efficient for computations because taking the trace of products of operators involving the Green's function $Q^{-1}$ are non-local. In \cite{P6} we presented a general method giving local representatives of the bivariant Chern character, by \emph{renormalization of the eta-cochain}. Let us recall how this can be used to calculate the diagonal of the commutative diagram (\ref{rrg}) in odd degree:
\be
\vcenter{\xymatrix{
\Kt_1(\Ic\hotimes\Ac) \ar[r]^{\rho_!} \ar[d] \ar@{.>}[rd]^{\Delta} & \Kt_1(\Ic\hotimes\Bc) \ar[d] \\
HP_1(\Ac) \ar[r]_{\ch(\rho)} & HP_1(\Bc) }}  \label{diag} 
\ee
A topological $K$-theory class $[u]\in \Kt_1(\Ic\hotimes\Ac)$ is represented by an invertible element $u\in (\Ic\hotimes\Ac)^+$ such that $u-1\in \Ic\hotimes\Ac$ (as usual $^+$ denotes unitalization of algebras). The diagonal of the diagram thus carries $[u]$ to the periodic Chern character $\ch(\rho_!(u))\in HP_1(\Bc)$. If one wants to use the non-local formulas, the diagonal is equivalently calculated by a cup-product between the Chern character $\ch(u)\in HP_1(\Ac)$ and the bivariant Chern character of the quasihomomorphism:
\be
\ch(\rho_!(u)) = \ch(\rho)\cdot\ch(u)\ .
\ee
In order to renormalize and get local formulas it will be necessary to impose that $u-1$ belongs to the dense subalgebra of finite size matrices $M_{\infty}(\Ac_0)\subset \Ic\hotimes\Ac$. Indeed, the algebra $\Ic\hotimes\Ac$ is too complete in general and the renormalization procedure used in the sequel works only with appropriate ``smooth'' elements. By \cite{CQ1}, the Chern character $\ch(u)\in HP_1(\Ac)$ is represented by the cycle of odd degree in the complex $X(\Th\Ac)$
\be
\ch_1(\uh)= \frac{1}{\sqrt{2\pi i}} \, \Tr\nat \uh^{-1}\dd\uh \quad \in \Om^1\Th\Ac_{\nat}\ ,
\ee
where $\uh$ is any lifting of $u$ to the unitalized (matrices over the) tensor algebra $M_{\infty}(\Th\Ac)^+$. If one chooses the canonical lifting $\uh=u$ induced by the linear inclusion $\Ac\hookrightarrow \Th\Ac$, then the image of $\ch_1(\uh)$ under the homotopy equivalence $\gamma: X(\Th\Ac)\to \Omh\Th\Ac$ is the $(b+B)$-cycle \cite{P5}
$$
\gamma \ch_1(\uh)= \frac{1}{\sqrt{2\pi i}} \sum_{n\geq 0}(-)^n n! \Tr(\uh^{-1}\dd\uh (\dd\uh^{-1}\dd\uh)^n)\ .
$$
For any choice of integer $n\geq 1$, the evaluation of this differential form under the chain map $\chih^{2n}\rho_*:\Omh\Th\Ac\to X(\Th\Bc)$ yields the representative $\ch^{2n}(\rho)\cdot\ch_1(\uh)$ of the class $\ch(\rho_!(u))$. Let us define $V= \rho_*\uh^{-1}[F,\rho_*\uh]$ as an element of the ideal $\Ic^s\hotimes\Th\Bc \subset \Mch^s$, and the Maurer-Cartan form $\om=\rho_*\uh^{-1}\dd(\rho_*\uh)$. Here we identify $\Lc(H_{\al})$ with $M_N(\cc)\otimes\Lc(H_{\al})$ for $N$ large enough in order to get rid of matrices. A straightforward computation gives
$$
\ch^{2n}(\rho)\cdot \ch_1(\uh) = \frac{1}{\sqrt{2\pi i}} \frac{(n!)^2}{(2n)!} \, \frac{1}{2}\tau\nat (V^{2n+2}\om-FV^{2n}\dd V) \quad \in \Om^1\Th\Bc_{\nat}
$$
for any $n\geq 1$. For $n>1$ one has $\frac{1}{2}\tau\nat (V^{2n+2}\om-FV^{2n}\dd V) = \tau\nat(V^{2n}\om)$, but this simplification does not hold in degree $n=1$ because the supertrace $\tau$ is defined only when $V$ is raised to a power $\geq 3$. Two consecutive degrees are related by the transgression relation 
\be
\ch^{2n}(\rho)\cdot \ch_1(\uh) - \ch^{2n+2}(\rho)\cdot \ch_1(\uh) = \nat\dd \big(\tch^{2n+1}(\rho)\cdot \ch_1(\uh)\big) \ . \label{tra}
\ee
Since the boundary map $\nat\dd: \Th\Bc \to \Om^1\Th\Bc_{\nat}$ factors through the commutator quotient space $\Th\Bc_{\nat}=\Th\Bc/[\ ,\ ]$ it is enough to compute $\tch^{2n+1}(\rho)\cdot \ch_1(\uh)$ modulo commutators and one finds 
$$
\nat\tch^{2n+1}(\rho)\cdot \ch_1(\uh) = \frac{1}{\sqrt{2\pi i}} \frac{(n!)^2}{(2n+1)!} \, \frac{1}{2}\tau\nat (FV^{2n+1})\quad \in \Th\Bc_{\nat}
$$
for $n\geq 1$. Of course the chain $\tch^{1}(\rho)\cdot \ch_1(\uh)$ corresponding to $n=0$ does not exist because $FV\in \Ic^s\hotimes\Th\Bc$ is not trace class. The trick used in \cite{P6} is to renormalize this quantity by any method, and write the Chern character as the boundary of a formal power series. Let us denote in the form of a $2\times 2$ matrix $\rho_*\uh= \bigl( \begin{smallmatrix} \uh_+ & 0 \\ 0 & Q^{-1}\uh_-Q \end{smallmatrix} \bigr) $, with $\uh_+$ and $Q^{-1}\uh_-Q$ in $\Lc(H_{\al})\hotimes\Th\Bc$. Then $V$ reads 
$$
V= \left( \begin{matrix} 0 & (1+Q^{-1}A)^{-1}-1 \\ Q^{-1}A & 0 \end{matrix} \right)\quad \mbox{with} \quad A:=\uh_-^{-1}Q \uh_+ - Q\ .
$$
One has $Q^{-1}A \in \Ic\hotimes\Th\Bc$ by Lemma \ref{lcom}. A simple calculation yields
$$
\nat\tch^{2n+1}(\rho)\cdot \ch_1(\uh) = \frac{(-)^n}{\sqrt{2\pi i}} \frac{(n!)^2}{(2n+1)!} \, \Tr\nat \left( \left(\frac{Q^{-1}A}{1+Q^{-1}A}\right)^{2n+1}(1+Q^{-1}A/2) \right)\ .
$$
$\Tr$ is the operator trace on $H_{\al}$. This quantity may be viewed as a formal power series in $A$ by taking the Neumann expansion $(1+Q^{-1}A)^{-1} = \sum_{k\geq 0}(-Q^{-1}A)^k$. Of course the series is divergent in general but it contains all the information we need. Writing the corresponding series for $n=0$, without taking care of traceability we would get
$$
\nat\tch^{1}(\rho)\cdot \ch_1(\uh) = \frac{1}{\sqrt{2\pi i}} \, \big( \Tr\nat(Q^{-1}A) - \frac{1}{2}\Tr\nat(Q^{-1}AQ^{-1}A) + \mbox{degree} \geq 3 \big)\ .
$$
The terms of degree $\geq 3$ are well-defined because the trace is finite. Thus only the first two terms need renormalization. This will be done explicitly in the next section, and this point uses the fact that $u-1\in M_{\infty}(\Ac_0)$. Once this renormalization is performed, we denote by $\nat\tch^{2n+1}_R(\rho)\cdot \ch_1(\uh)$ the formal power series obtained for any values of $n\geq 0$. The first term of this series has degree $2n+1$ in $A$, and thus increases with $n$. As a consequence, the infinite sum
\be
W_R(A) := \sqrt{2\pi i} \, \sum_{n\geq 0} \nat\tch^{2n+1}_R(\rho)\cdot \ch_1(\uh) \label{W}
\ee
makes perfectly sense as a formal power series in $A$. It is not convergent in general. Proceeding as in \cite{P5} Lemma 4.2, one finds that the term of degree $k$ in this expansion is simply $W^k(A)= \frac{(-)^{k+1}}{k} \Tr\nat((Q^{-1}A)^k)$ for $k\geq 3$, and in lower degrees 
$$
W_R^1(A) = \Tr\nat(Q^{-1}A) \quad \mbox{and} \quad W_R^2(A)=-\frac{1}{2} \Tr\nat(Q^{-1}AQ^{-1}A)
$$
are renormalized quantities. Note that renormalization is not unique, and in general $W_R(A)$ is defined modulo addition of a ``local'' polynomial $P(A)$ of degree 2. The series $W_R(A)$ can be depicted in terms of Feynman graphs as follows. We represent the Green's operator by an arrow $Q^{-1}=\xymatrix{*{}\ar@{-}[r] |-{\SelectTips{cm}{}\object@{>}}  & *{}}$ and the insertion of potential $A=\bullet$ by a point. The products $Q^{-1}AQ^{-1}A\ldots$ are represented by chains, which are closed by taking the trace $\Tr\nat$. Hence $W_R(A)$ is given by a formal series of loops \\
$$
W_R(A)  =\ \xymatrix{ *=0{\bullet} \ar@(ur,dr)@{-}[] |-{\SelectTips{cm}{}\object@{>}} }   
\ - \frac{1}{2}\  
\vcenter{\xymatrix@R=3pc{ *=0{\bullet} \ar@/^/@{-}[d] |-{\SelectTips{cm}{}\object@{>}} \\*=0{\bullet} \ar@/^/@{-}[u] |-{\SelectTips{cm}{}\object@{>}} } }
\ +\ \frac{1}{3}\ 
\vcenter{\xymatrix@C=1pc{ & *=0{\bullet} \ar@{-}[dr] |-{\SelectTips{cm}{}\object@{>}} &  \\
*=0{\bullet} \ar@{-}[ur] |-{\SelectTips{cm}{}\object@{>}} & & *=0{\bullet} \ar@{-}[ll] |-{\SelectTips{cm}{}\object@{>}} }}
\ -\frac{1}{4}\ 
\vcenter{\xymatrix@R=3pc@C=2.8pc{ *=0{\bullet} \ar@{-}[r] |-{\SelectTips{cm}{}\object@{>}} & *=0{\bullet} \ar@{-}[d] |-{\SelectTips{cm}{}\object@{>}} \\
*=0{\bullet} \ar@{-}[u] |-{\SelectTips{cm}{}\object@{>}} & *=0{\bullet} \ar@{-}[l] |-{\SelectTips{cm}{}\object@{>}} }}
\ +\frac{1}{5}\ 
\vcenter{\xymatrix@R=0.25pc@C=0.5pc{ & & *=0{\bullet} \ar@{-}[ddrr] |-{\SelectTips{cm}{}\object@{>}} & \\
 & & & & \\
*=0{\bullet} \ar@{-}[uurr] |-{\SelectTips{cm}{}\object@{>}} &  &  &  & *=0{\bullet} \ar@{-}[dddl] |-{\SelectTips{cm}{}\object@{>}} \\
 & & & &  \\
 & & & & \\
 & *=0{\bullet} \ar@{-}[uuul] |-{\SelectTips{cm}{}\object@{>}} & & *=0{\bullet} \ar@{-}[ll] |-{\SelectTips{cm}{}\object@{>}} & }}
\ +\ldots
$$
Each term of this series lies in $\Th\Bc_{\nat}$. Its universal derivative under the boundary map $\dd:\Th\Bc_{\nat}\to \Om^1\Th\Bc_{\nat}$ may be computed by means of the BRS transformations \cite{MSZ}. We make the convention that the operators $A$ and $Q$ have odd parity, which means that $\dd$ anticommutes with $Q$, and the derivative of $A=\uh_-^{-1}Q \uh_+ - Q$ reads
\be
\dd A= -(Q+A)\om_+-\om_-(Q+A)\ ,
\ee
with the $2\times 2$ matrix notation $\om= \bigl( \begin{smallmatrix} \om_+ & 0 \\ 0 & Q^{-1}\om_-Q \end{smallmatrix} \bigr)$. A simple computation (see \cite{P6}) shows that in the sense of formal power series the derivatives of graphs of order $k \geq 3$ cancel each other, whereas the low degree graphs yield an \emph{anomaly} $\Delta(\om,A)$ polynomial in $A$ and linear in $\om$. The anomaly is necessarily a cycle in $\Om^1\Th\Bc_{\nat}$ and thus defines a class of odd degree $\Delta(\om,A)\in HP_1(\Bc)$. Diagramatically: 
$$
\vcenter{\xymatrix@!0@=2.5pc{
W_R(A) \ar[d]_{\dd} & = &  & W^1_R(A) \ar[dl] \ar[dr]  & + & W^2_R(A) \ar[dl] \ar[dr] & + & W^3(A) \ar[dl] \ar[dr] & + & W^4(A) \ar[dl] \ar[dr] & +  \ldots  \\
\Delta(\om,A) & = & \Delta^0(\om,A) & + & \Delta^1(\om,A) & + & \Delta^2(\om,A) & + & 0 & + & \ 0  \ldots }}
$$
Each component of the anomaly $\Delta^k(\om,A)$ is a homogeneous polynomial of degree $k$ in $A$. At the same time, taking the boundary of the series (\ref{W}) and using the transgression relation (\ref{tra}) valid for $n\geq 1$ yields  
$$
\dd W_R(A) = \sqrt{2\pi i}\, \big( \nat\dd \big(\tch_R^1(\rho)\cdot \ch_1(\uh)\big) + \ch^2(\rho)\cdot \ch_1(\uh) \big)
$$
and shows that the anomaly $\Delta(\om,A)$ defines the same cyclic homology class as $\sqrt{2\pi i}\, \ch(\rho)\cdot \ch(u)$ in $HP_1(\Bc)$. Thus we have proved

\begin{proposition}\label{pano}
Let $u\in M_{\infty}(\Ac_0)^+$ be an invertible element representing a class in $\Kt_1(\Ac)$. Denote by $\rho_*\uh= \bigl( \begin{smallmatrix} \uh_+ & 0 \\ 0 & Q^{-1}\uh_-Q \end{smallmatrix} \bigr) $ the image of its canonical lifting $\uh$ under the homomorphism $\rho_*:\Th\Ac\to\Mch^s_+$. For any renormalization of the low-order graphs of the formal series $W_R(A)$ associated to the ``gauge potential'' $A=\uh_-^{-1}Q \uh_+ - Q$, the anomaly
\be
\Delta(\om,A) \equiv \sqrt{2\pi i}\, \ch(\rho_!(u)) \ \in HP_1(\Bc)
\ee
is a polynomial of degree at most 2 in $A$ which computes the diagonal of the commutative diagram (\ref{diag}), up to a factor $\sqrt{2\pi i}$. \cqfd
\end{proposition}

Let us write explicitly the elements $A$ and $\om$ in terms of the given invertible $u\in M_{\infty}(\Ac_0)^+$. By hypothesis, the difference $u-1$ lies in $M_{\infty}(\Ac_0)$ hence it is a finite sum 
$$
u-1 = \sum_{g\in G} u(g)U^*_g\ ,
$$
with components $u(g)$ in the matrix algebra over $\cinfc(\Si)$. Since $u(g)$ is represented by a bounded operator on $H_{\al}$ by pointwise multiplication, the canonical lifting $\uh_+=\rho_*(\uh)_+$ reads
$$
\uh_+-1= \sum_{g\in G} u(g)r(g)_+ \otimes U^*_g\ \in \Lc(H_{\al})\hotimes \Th\Bc\ ,
$$
where $r(g)_+$ is the representation of $G$ by pullback. Recall also that we have a ``minus'' representation $r(g)_-$ acting on the space of one-forms over $\Si$, and the lifting $\uh_-=Q \rho_*(\uh)_- Q^{-1}$ is given by an analogous formula. We may rewrite this as
$$
\uh_{\pm}-1 = \sum_{g\in G}\uh(g)r(g)_{\pm} \quad\mbox{with}\quad \uh(g):=u(g)\otimes U^*_g \in M_{\infty}(\cinfc(\Si))\otimes\Th\Bc\ ,
$$
and the algebraic tensor product $\cinfc(\Si)\otimes\Th\Bc  = \varprojlim_k \cinfc(\Si)\otimes T\Bc/(J\Bc)^k$ is defined in the sense of pro-algebras. The explicit formulas for the inverses $\uh^{-1}_{\pm}$ are more complicated. Setting $u^{-1}-1=\sum_{g\in G} u^{-1}(g) U^*_g$ one finds a series
\beq
\lefteqn{ \uh^{-1}_{\pm} - 1 = \sum_{g\in G} u^{-1}(g)r(g)_{\pm} \otimes  U^*_g + } \non\\
&& \mspace{-40mu} \sum_{n\geq 1} \Big(1+ \sum_{g\in G} u^{-1}(g)r(g)_{\pm} \otimes  U^*_g\Big) \Big( \sum_{h,i\in G} u(h)u^{-1}(i)^hr(ih)_{\pm} \otimes ( U^*_{ih} - U^*_h \otimes  U^*_i) \Big)^n \non
\eeq
where $u^{-1}(i)^h$ is the pullback of the function $u^{-1}(i)\in M_{\infty}(\cinfc(\Si))$ by the conformal mapping $h$. Since $U^*_{ih} - U^*_h \otimes  U^*_i$ is in the ideal $J\Bc$ the sum over $n$ converges for the $J\Bc$-adic topology of $\Th\Bc$. As before we can decompose
$$
\uh_{\pm}^{-1}-1 = \sum_{g\in G}\uh^{-1}(g)r(g)_{\pm} \ ,\qquad \uh^{-1}(g) \in M_{\infty}(\cinfc(\Si))\otimes\Th\Bc\ ,
$$
but now the components $\uh^{-1}(g)$ are given by complicated expressions and the sum over $g\in G$ is pro-finite. Let us compute, taking the equivariance relation $Qr(g)_+=r(g)_-Q$ into account, 
$$
Q\uh_+-\uh_-Q = \sum_{g\in G}(Q\uh(g)r(g)_+ - \uh(g)r(g)_-Q)=\sum_{g\in G} Q\cdot\uh(g)\, r(g)_+\ .
$$
One has $Q\cdot\uh(g)= (Q\cdot u(g))\otimes U^*_g \in M_{\infty}(\Omc^{0,1}(\Si))\otimes \Th\Bc$ and $Q\cdot u(g)=d\zb\d_{\zb} u(g)$ is a matrix-valued one-form viewed as an operator. This allows to write explicitly the gauge potential $A=\uh_-^{-1}Q\uh_+ - Q$ as
$$
A=\sum_{g\in G} A(g) r(g)_+\ ,\qquad A(g)=\sum_{h\in G}\uh^{-1}(h)(Q\cdot\uh(gh^{-1}))^h \in M_{\infty}(\Omc^{0,1}(\Si))\otimes\Th\Bc
$$
where $^h$ denotes the pullback of differential forms by the conformal mapping $h$. The fact that $A(g)$ is a (matrix of) smooth compactly supported differential form with values in $\Th\Bc$ enables to view the $k$-linear functionals $W^k(A)$ as \emph{distributions}. Finally, the Maurer-Cartan form may also be decomposed,
$$
\om_{\pm}=\sum_{g\in G} \om(g) r(g)_{\pm}\ ,\quad \om(g)=\sum_{h\in G}\uh^{-1}(h)(\dd\uh(gh^{-1}))^h \in M_{\infty}(\cinfc(\Si))\otimes \Om^1\Th\Bc
$$
and the BRS transformation $\dd A= -(Q+A)\om_+-\om_-(Q+A)$ explicitly reads
\be
\dd A = -\sum_{g\in G} Q\cdot\om(g)\, r(g)_+ - \sum_{g,h\in G} \big(A(h)\om(g)^h + \om(h)A(g)^h\big) r(gh)_+ \ . \label{BRS}
\ee

\section{Conformal renormalization}\label{sren}

In this section we renormalize the low-order terms of the formal power series $W_R(A)$ associated to a given gauge potential 
$$
A= \sum_{g\in G} A(g)r(g)_+\ ,\qquad A(g)\in M_{\infty}(\Omc^{0,1}(\Si))\otimes\Th\Bc\ ,
$$
and calculate the corresponding anomaly. We write $A(g)=d\zb A_{\zb}(g)$ and view the component $A_{\zb}(g)\in M_{\infty}(\cinfc(\Si))\otimes \Th\Bc$ as a matrix of smooth compactly supported functions with values in the pro-algebra $\Th\Bc$. As before we denote by $\Tr$ the operator trace on $\Lc(H_{\al})$ and by $\nat: \Th\Bc\to\Th\Bc_{\nat}$ the projection onto the commutator quotient space. The $k$-th term of the series 
$$
W^k(A)=\frac{(-)^{k+1}}{k}\Tr\nat(Q^{-1}A)^k\ \in \Th\Bc_{\nat}
$$ 
should be defined as a $k$-linear functional of the test functions $A_{\zb}(g)$, i.e. a distribution, with values in $\Th\Bc_{\nat}$. For $k\geq 3$ we know that the operator trace is finite and $W^k(A)$ is indeed a distribution. For $k=1,2$ the task is therefore to give a distributional meaning to the ill-defined quantities
$$
W_R^1(A)= \Tr\nat(Q^{-1}A) \ ,\qquad W_R^2(A)= -\frac{1}{2} \Tr\nat(Q^{-1}AQ^{-1}A) \ .
$$
The operators $Q^{-1}$ and $A(g)r(g)_+$ have distribution kernels whose evaluation at two complex points $(z,w)\in \Si\times \Si$ reads
$$
Q^{-1}(z,w)=\frac{1}{\pi(z-w)}\ ,\qquad \big(A(g)r(g)_+\big)(z,w)= A_{\zb}(g,z) \delta^2(g(z)-w)\ .
$$
$A_{\zb}(g,z) \in \Th\Bc$ denotes the evaluation of the function $A_{\zb}(g)$ on point $z\in \Si$, and $\delta^2(g(z)-w)$ is the Dirac current associated to the submanifold $w=g(z)$ of $\Si\times \Si$. Also recall that by definition $g$ is a holomorphic function of $z$, and $\supp A(g) \in \dom(g)$. \\

The renormalization procedure will be based on the following fundamental result. Fix a point $z_0\in \Si$ and for some integer $n\geq 1$ consider the meromorphic function 
$$
f(z) = \frac{1}{(z-z_0)^n}\ .
$$
Then $f$ is smooth over $\Si\backslash \{z_0\}$ and has an isolated singularity at $z=z_0$. Our task is to extend $f$ to a distribution over the entire plane $\Si$. This can be done as follows. For $n=1$, the Cauchy kernel $1/(z-z_0)$ is a locally integrable function over $\Si$, hence it canonically defines a distribution. According to the terminology of Epstein and Glaser \cite{EG}, its singularity order at $z_0$ is $-1$. There is no ambiguity in that case. For $n=2$ unfortunately, the function $1/(z-z_0)^2$ is not a distribution. However one has, for $z\neq z_0$, the equality of functions
$$
\frac{1}{(z-z_0)^2} = -\d_z \left(\frac{1}{z-z_0}\right) \ ,
$$
and the right-hand-side makes sense as a distribution over the entire plane. Thus we can \emph{define} the distributional extension of $f$ accordingly. Note that its singularity order is $0$. Now this extension is not unique. Indeed if $\ft$ denotes another extension, the difference $\ft-f$ is a distribution with support localized at $\{z_0\}$, hence it is necessarily proportional to a finite sum of derivatives of the Dirac measure $\delta^2(z-z_0)$. If one requires to preserve the singularity order, one thus finds that all the possible distributional extensions are given by 
\be
\frac{1}{(z-z_0)^2} \rightsquigarrow -\d_z \frac{1}{z-z_0} + a\delta^2(z-z_0) \label{e1}
\ee
for some parameter $a$. Among all these possibilities only the choice $a=0$ is \emph{conformally invariant}. Indeed if $w=h(z)$ is any other complex coordinate system around the singularity, with $h$ a holomorphic function and $w_0=h(z_0)$, one has the equality of distributions
$$
\d_z \left( \frac{1}{z-z_0} \right) \left(\frac{z-z_0}{w-w_0}\right)^2 = \d_w \left( \frac{1}{w-w_0} \right)\ .
$$
This shows that for $a=0$ the extension only depends on the complex structure of $\Si$. The situation for $n>2$ is analogous: all the possible distributional extensions of $1/(z-z_0)^n$ with singularity order $n-2$ are given by
\be
\frac{1}{(z-z_0)^n} \rightsquigarrow \frac{(-)^{n-1}}{(n-1)!} \d_z^{n-1} \left(\frac{1}{z-z_0}\right) + P(\d_z,\d_{\zb})\delta^2(z-z_0)\ , \label{e2}
\ee
where $P(\d_z,\d_{\zb})$ is a polynomial of degree $n-2$ in the derivatives. Only the choice $P=0$ yields a conformally invariant distribution: for any complex coordinate system $w$ one has the equality of distributions
$$
\d_z^{n-1} \left( \frac{1}{z-z_0} \right) \left(\frac{z-z_0}{w-w_0}\right)^n = \d_w^{n-1} \left( \frac{1}{w-w_0} \right)\ .
$$
This kind of renormalization scheme is well-known in conformal field theory models. In our case the simple formulas (\ref{e1}) and (\ref{e2}) are sufficient to renormalize the formal power series $W_R(A)$. The general theory of renormalization by distributional extension (or distributional splitting in the Minkowskian framework) is due to Epstein and Glaser \cite{EG}.\\

Let us start with the tadpole graph $W_R^1(A)$. The operator trace of $Q^{-1}A$ is obtained by integrating its distribution kernel along the diagonal. Disregarding traceability, one gets formally
$$
W_R^1(A) = \sum_{g\in G} \int\!\!\!\!\int d^2z\,d^2w\, \frac{\nat A_{\zb}(g,z)}{\pi(w-z)}\, \delta^2(g(z)-w) = \sum_{g\in G}\int d^2z\, \frac{\nat A_{\zb}(g,z)}{\pi(g(z)-z)}\ .
$$
This expression will make sense once we interpret the function $z\mapsto 1/(g(z)-z)$ as a distribution. We see that the ambiguity comes from the singular points where the function $g(z)-z$ vanishes, i.e. at the fixed point set of the mapping $g$. Renormalizing $W_R^1(A)$ thus amounts to extend the function $1/(g(z)-z)$ to a distribution in the neighborhood of all fixed points. The Taylor expansion of the holomorphic function $g$ around a given fixed point $z_0=g(z_0) \in\dom(g)$ leads to
$$
g(z)-z= (z-z_0)(g'(z_0)-1)+ \frac{(z-z_0)^2}{2}g''(z_0) + \ldots + \frac{(z-z_0)^n}{n!}g^{(n)}(z_0) + \ldots
$$
We say that the fixed point $z_0$ is of order $n$ (with $1\leq n\leq \infty$) if $g(z)-z$ behaves like $(z-z_0)^n$ when $z\to z_0$, that is, if the first non-vanishing term in the Taylor expansion is $(z-z_0)^n$. Hence a fixed point of order 1 verifies $g'(z_0)\neq 1$; in other words its tangent space is not fixed by the tangent map. A fixed point of order $\geq 2$ has a fixed tangent space ($g'(z_0)=1$), and $g$ differs from the identity only by the higher order jets. A fixed point has order $n=\infty$ if and only if $g$ is the identity mapping in a neighborhood of $z_0$. Note that for a fixed point $z_0$ of order $n< \infty$, the function
\be
H^n_{g,z_0}(z) := \frac{(z-z_0)^n}{g(z)-z} \label{fun}
\ee
is holomorphic in a neighborhood of $z_0$. Now the distributional extension of the function $1/(g(z)-z)$ depends on the order $n$. Three distinct cases can occur:\\

\noindent {\bf i) $n=1$:} Then $z_0$ is necessarily an isolated fixed point and one has the equality of functions
$$
\frac{1}{g(z)-z} = \frac{1}{z-z_0} \, H^1_{g,z_0}(z)
$$
for $z\neq z_0$. We know that $1/(z-z_0)$ is already a well-defined distribution on a neighborhood of $z_0$, hence also $1/(g(z)-z)$ and the quantity $W_R^1(A)$ does not need any renormalization in this case.\\

\noindent {\bf ii) $2\leq n < \infty$:} Then $z_0$ is necessarily an isolated fixed point and one has the equality of functions
$$
\frac{1}{g(z)-z} = \frac{1}{(z-z_0)^n} \, H^n_{g,z_0}(z)
$$
for $z\neq z_0$. According to the discussion above we perform the distributional extension
\be
\frac{1}{g(z)-z}\ \rightsquigarrow\ \frac{(-)^{n-1}}{(n-1)!}\, \d_z^{n-1} \left(\frac{1}{z-z_0}\right) \, H^n_{g,z_0}(z)\ . \label{ext1}
\ee
It is not unique and one could also add a linear combination of the Dirac measure $\delta^2(z-z_0)$ and its derivatives up to order $n-2$. \\

\noindent {\bf iii) $n=\infty$:} Then $g(z)=z$ in a neighborhood of $z_0$, hence all the points are fixed. In this case the function $1/(g(z)-z)$ has no meaning at all, and renormalization consists in replacing the function by any quantity, for example zero:
\be
\frac{1}{g(z)-z}\ \rightsquigarrow\ 0\ .
\ee
Again this choice is non-unique and any other numerical value would work as well.\\

Remark that the distributional extension given in case ii) also gives the correct answer when $n=1$. We deal now with the second graph $W_R^2(A)$. Integrating the distribution kernel of the operator $Q^{-1}AQ^{-1}A$ along the diagonal one formally gets 
\beq
W_R^2(A) &=& -\frac{1}{2} \sum_{g,h\in G}   \int d^2v\, d^2z\, d^2w\,\nat \frac{A_{\zb}(g,z)}{\pi(v-z)}\, \frac{A_{\wb}(h,w)}{\pi(g(z)-w)}\, \delta^2(h(w)-v) \non\\
&=& -\frac{1}{2} \sum_{g,h\in G} \int d^2z\, d^2w\, \frac{\nat A_{\zb}(g,z) A_{\wb}(h,w)}{\pi^2(h(w)-z)(g(z)-w)}  \ . \non
\eeq
Again, this expression will make sense once we extend the function of two variables $(z,w)\mapsto 1/((h(w)-z)(g(z)-w))$ to a distribution in the neighborhood of singular points. Since $h$ is a diffeomorphism, the function $1/(h(w)-z)$ actually defines unambiguously a distribution of two variables, with singular support restricted to the diagonal $h(w)=z$. The same is true for $1/(g(z)-w)$. Hence only the product of these two distributions may be ill-defined on their common singular support, i.e. at points $(z_0,w_0)$ such that $h(w_0)=z_0$ and $g(z_0)=w_0$. This implies that $z_0=hg(z_0)$ and $w_0=gh(w_0)$ must be fixed points, and the task is to extend the function $1/((h(w)-z)(g(z)-w))$ to a distribution in the neighborhood of $(z_0,w_0)$. As before this depends on the order $n$ of the fixed point $z_0$ (or equivalently $w_0$):\\

\noindent {\bf iv) $n=1$:} Then $z_0$ is an isolated fixed point of $hg$, and $w_0$ an isolated fixed point of $gh$. We perform the change of variable $w=g(v)$ around $w_0$ and write
\beq
\frac{1}{(h(w)-z)(g(z)-w)} &=& \frac{1}{(hg(v)-z)(g(z)-g(v))} \non\\
&=& \frac{1}{(hg(v)-z)(z-v)} \left( \frac{z-v}{g(z)-g(v)} \right) \non\\
&=& \frac{1}{(hg(v)-v)}  \left( \frac{1}{hg(v)-z}+\frac{1}{z-v} \right)\frac{z-v}{g(z)-g(v)} \non\ .
\eeq
Since $g$ is a diffeomorphism, the quotient $(z-v)/(g(z)-g(v))$ is a holomorphic function of the complex variable $z-v$. By hypothesis the fixed point $v_0=z_0$ is of order 1, hence the function $hg(v)-v$ is equivalent to $(v-v_0)((hg)'(v_0)-1)$ when $v\to v_0$, and the inverse function $1/(hg(v)-v)$ unambiguously extends to a distribution in a neighborhood of $v_0$. Moreover, the complex variables $hg(v)-v$ and $z-v$ are independent, so the product of distributions $1/(hg(v)-v)$ and $1/(z-v)$ is well-defined. Similarly for the product of $1/(hg(v)-v)$ and $1/(hg(v)-z)$. \\

\noindent {\bf v) $2\leq n < \infty$:} Then $z_0$ and $w_0$ are isolated, and by the same change of variable $w=g(v)$ one still has
$$
\frac{1}{(h(w)-z)(g(z)-w)} = \frac{1}{(hg(v)-v)}  \left( \frac{1}{hg(v)-z}+\frac{1}{z-v} \right) \frac{z-v}{g(z)-g(v)} \non\ .
$$
Now the function $1/(hg(v)-v)$ does not extend canonically to a distribution at the fixed point $v_0=z_0$. Proceeding as before we renormalize
\be
\frac{1}{(hg(v)-v)}\ \rightsquigarrow\ \frac{(-)^{n-1}}{(n-1)!}\, \d_v^{n-1} \left(\frac{1}{v-v_0}\right) \, H^n_{hg,v_0}(v)\ ,
\ee
and the r.h.s. makes sense as a distribution in the variable $v$. For the same reason as case iv) taking products with $1/(z-v)$ or $1/(hg(v)-z)$ also yields well-defined distributions. It is interesting to investigate the freedom one has in the choice of this extension. We know that the distributional extension of $1/(hg(v)-v)$ is not unique and we can add a linear combination of $\delta^2(v-v_0)$ and its derivatives up to order $n-2$. But these local terms will not contribute to the final result because 
$$
\d^k_v \d^l_{\vb} \delta^2(v-v_0)\, \left( \frac{1}{hg(v)-z}+\frac{1}{z-v} \right) = 0
$$
in the sense of distributions whenever $k+l\leq n-1$ (indeed $(hg)'(v_0)=1$ and the higher derivatives of $hg$ vanish at $v_0$ up to order $n-1$). Hence unlike the tadpole graph situation, the renormalization of $W_R^2(A)$ is unique in the case of isolated fixed points. \\

\noindent {\bf vi) $n=\infty$:} It means that $hg(z)=z$ in a neighborhood of $z_0$, or equivalently $gh(w)=w$ in a neighborhood of $w_0$. As usual use the change of variables $w=g(v)$: 
$$
\frac{1}{(h(w)-z)(g(z)-w)} = \frac{1}{(v-z)(g(z)-g(v))} = \frac{-1}{(z-v)^2} \left( \frac{z-v}{g(z)-g(v)} \right) \ .
$$
The quotient $(z-v)/(g(z)-g(v))$ is a holomorphic function of $z-v$, but the function $1/(z-v)^2$ does not extend canonically to a distribution on the diagonal $z=v$. We renormalize
\be
\frac{-1}{(z-v)^2}\ \rightsquigarrow\ \d_z \frac{1}{z-v}\ ,
\ee
so that everything makes sense. In this case the renormalization is not unique, it is always possible to add a term proportional to the Dirac measure $\delta^2(z-v)$. \\

The higher order graphs $W^k(A)$ for $k\geq 3$ are unambiguously defined because the operator trace converges. For example the triangle graph may be evaluated as
$$
W^3(A) = \frac{1}{3} \sum_{g,h,i\in G} \int d^2z\, d^2w\, d^2v\, \frac{\nat A_{\zb}(g,z) A_{\wb}(h,w) A_{\vb}(i,v)}{\pi^3(i(v)-z)(g(z)-w)(h(w)-v)} \ ,
$$
and the function $(z,w,v)\mapsto 1/(i(v)-z)(g(z)-w)(h(w)-v)$ extends to a unique distribution with the same singularity order on the diagonal. This can be seen as follows. Perform the change of variables $w=g(s)$, $v=hg(u)$ and write outside the diagonal
\beq
\lefteqn{\frac{1}{(i(v)-z)(g(z)-w)(h(w)-v)} = \frac{1}{(ihg(u)-z)(g(z)-g(s))(hg(s)-hg(u))} } \non\\
&=& \frac{1}{(ihg(u)-z)(z-s)}\, \frac{z-s}{g(z)-g(s)}\, \frac{1}{hg(s)-hg(u)} \non\\
&=& \left(\frac{1}{ihg(u)-z}+\frac{1}{z-s}\right) \frac{1}{(ihg(u)-s)(s-u)}\, \frac{s-u}{hg(s)-hg(u)}\, \frac{z-s}{g(z)-g(s)} \ .\non 
\eeq
$(s-u)/(hg(s)-hg(u))$ and $(z-s)/(g(z)-g(s))$ are holomorphic functions. $1/(ihg(u)-z)$ and $1/(z-s)$ are well-defined distributions. It is therefore sufficient to find a distributional extension to the function $1/((ihg(u)-s)(s-u))$. This was already done in items iv) - vi) above, and depends on the nature of the fixed points $u_0$ of the mapping $ihg$. If $u_0$ is isolated then we know that the extension is unique. Otherwise if $ihg(u)=u$ in the neighborhood of $u_0$ then all possible extensions (with the same singularity degree) are of the form
$$
\frac{-1}{(u-s)^2}\ \rightsquigarrow\ \d_u \frac{1}{u-s} + a\delta^2(u-s)\ ,
$$
for some parameter $a$. But the ambiguity carried by the Dirac measure does not affect the final result, because
$$
\left(\frac{1}{u-z}+\frac{1}{z-s}\right) \delta^2(u-s)=0\ .
$$
Hence $W^3(A)$ is uniquely defined as a distribution.\\

Let us now calculate the anomaly, i.e. the image of the formal power series
$$
W_R(A) = W_R^1(A) + W_R^2(A) + W^3(A) + \ldots
$$
under the boundary map $\dd : \Th\Bc_{\nat} \to \Om^1\Th\Bc_{\nat}$. It amounts to derive each functional $W^k_R(A)$ according to the BRS transformation $\dd A=-(Q+A)\om_+ - \om_-(Q+A)$. In terms of the components $A(g) \in M_{\infty}(\Omc^{0,1}(\Si))\otimes \Th\Bc$ and $\om(g)\in M_{\infty}(\cinfc(\Si))\otimes \Om^1\Th\Bc$ Equation (\ref{BRS}) yields:  
\beq
\dd A_{\zb}(g,z) &=& \d_{\zb} \om(g,z) + \non\\
&& \sum_{h\in G} \big( A_{\zb}(h,z) \om(gh^{-1},h(z)) - \om(h,z) \d_{\zb}\overline{h(z)} A_{\hb}(gh^{-1},h(z))\big)  \non
\eeq
at any point $z\in \Si$. From Proposition \ref{pano} we know that the anomaly $\Delta(\om,A)=\Delta^0(\om,A)+\Delta^1(\om,A)+\Delta^2(\om,A)$ is a polynomial of degree at most 2 in $A$. Because the ambiguities of renormalization are localized at the fixed points of the conformal mappings, the anomaly itself is necessarily given by a formula localized at fixed points. 

\begin{proposition}\label{ploc}
For the conformal renormalization chosen above the anomaly is a polynomial of degree one in $A$. Its component of degree zero $\Delta^0(\om,A)$ is a sum over the isolated fixed points for all mappings $g\in G$:
\be
\Delta^0(\om,A)=\sum_{g\in G}\sum_{\substack{z_0=g(z_0)\\ \textup{isolated}}} \frac{-1}{(n-1)!}\, \d_z^{n-1}\big(H_{g,z_0}^n(z) \nat\om(g,z)\big)_{z=z_0}\ , \label{d0}
\ee
where $n$ denotes the order of the fixed point $z_0$. The component of degree one $\Delta^1(\om,A)$ is an integral over the manifold of non-isolated fixed points:
\be
\Delta^1(\om,A)=\frac{1}{\pi} \sum_{g,h\in G}\int_{z=hg(z)}d^2z\, \nat \big(\d_z- \frac{1}{2}\d_z\ln g'(z)\big)A_{\zb}(g,z)\,\,\om(h,g(z))\ . \label{d1}
\ee
\end{proposition}
{\it Proof:} Let us calculate the variation of the tadpole graph $W^1_R(A)$ under the BRS transformation. It is the sum of a degree zero term and a degree one term with respect to $A$,
\beq
\lefteqn{ \dd W_R^1(A) = \sum_{g\in G} \int d^2z\, \frac{1}{\pi(g(z)-z)} \, \nat\d_{\zb}\om(g,z)\ + } \non\\
&& \sum_{g,h\in G} \int d^2z\, \frac{1}{\pi(gh(z)-z)}\, \nat \big( A_{\zb}(h,z) \om(g,h(z)) - \om(h,z) \d_{\zb}\overline{h(z)} A_{\hb}(g,h(z))\big)  \non
\eeq
and the function $1/(g(z)-z)$ is extended to a distribution at singular points by i) - iii). The degree zero part of the anomaly comes entirely from the first term. Hence
$$
\Delta^0(\om,A)= \sum_{g\in G} \int d^2z\, \frac{1}{\pi(g(z)-z)} \, \nat\d_{\zb}\om(g,z)  \ ,
$$
and we can integrate by parts because this integral is the evaluation of a distribution on the test function $\d_{\zb}\om(g,z)$. Since $1/(g(z)-z)$ is a holomorphic function when $g(z)\neq z$, its $\d_{\zb}$ derivative vanishes but the singular points may contribute. Hence we can suppose that the support of $\om(g)$ is concentrated in the neighborhood of an isolated fixed point $z_0$ of order $n<\infty$ and calculate
\beq
\lefteqn{ -\int d^2z\, \d_{\zb} \left( \frac{(-)^{n-1}}{(n-1)!}\, \d_z^{n-1}\frac{1}{\pi(z-z_0)} \, H^n_{g,z_0}(z) \right) \, \nat\om(g,z) } \non\\
&&\qquad\qquad = \int d^2z\,  \frac{-1}{(n-1)!}\, \d_{\zb}\frac{1}{\pi(z-z_0)} \, \d_z^{n-1} \big( H^n_{g,z_0}(z) \, \nat\om(g,z) \big) \non\\
&&\qquad\qquad = \int d^2z\,  \frac{-1}{(n-1)!}\, \delta^2(z-z_0) \, \d_z^{n-1} \big( H^n_{g,z_0}(z) \, \nat\om(g,z) \big) \non\\
&&\qquad\qquad  = \frac{-1}{(n-1)!} \, \d_z^{n-1} \big( H^n_{g,z_0}(z) \, \nat\om(g,z) \big)_{z=z_0} \ . \non
\eeq
There is no contribution from the non-isolated fixed points because $1/(g(z)-z)$ is renormalized to zero in that case. Hence summing over all $g\in G$ and all isolated fixed points one gets
$$
\Delta^0(\om,A)=\sum_{g\in G}\sum_{\substack{z_0=g(z_0)\\ \textup{isolated}}} \frac{-1}{(n-1)!}\, \d_z^{n-1}\big(H_{g,z_0}^n(z) \nat\om(g,z)\big)_{z=z_0}\ .
$$
The second term of $\dd W_R^1(A)$, linear in $A$, will contribute to the degree one of the anomaly $\Delta^1(\om,A)$, but there is also another contribution from the BRS variation of $W_R^2(A)$. The latter is a sum of a linear term and a quadratic term in $A$:
\beq
\lefteqn{ \dd W_R^2(A) =  -\frac{1}{2\pi^2}\sum_{g,h\in G}\int d^2z\, d^2w\, \nat\big( \d_{\zb}\om(g,z) A_{\wb}(h,w) + A_{\zb}(g,z) \d_{\wb}\om(h,w) } \non\\
&& + [A,\om]_{\zb}(g,z) A_{\wb}(h,w) + A_{\zb}(g,z) [A,\om]_{\wb}(h,w) \big)/((h(w)-z)(g(z)-w))\ . \non
\eeq
One could be tempted to use the apparent symmetry $(g,z) \leftrightarrow (h,w)$ to simplify the numerator. However it is not obvious that the choice of distributional extension we made for the function $1/((h(w)-z)(g(z)-w))$ actually possesses this symmetry. Let us integrate by parts the contribution of $\dd W_R^2(A)$ to $\Delta^1(\om,A)$. As before only the singular points are relevant; hence we first investigate the role of an isolated fixed point $z_0=hg(z_0)$ of order $n<\infty$. The renormalization performed in iv) - v) yields 
$$
\frac{1}{(h(w)-z)(g(z)-w)} = \frac{1}{(hg(v)-v)}  \left( \frac{1}{hg(v)-z}+\frac{1}{z-v} \right) \frac{z-v}{g(z)-g(v)} 
$$
with the change of variables $w=g(v)$ and the distributional extension
$$
\frac{1}{(hg(v)-v)} = \frac{(-)^{n-1}}{(n-1)!}\, \d_v^{n-1} \left(\frac{1}{v-v_0}\right) \, H^n_{hg,v_0}(v) 
$$
around the fixed point $v_0=z_0$. Then we integrate by parts, keeping the relation $w=g(v)$ in mind:
$$
-\frac{1}{2}\int d^2z\, d^2w\, \frac{1}{\pi^2(h(w)-z)(g(z)-w)}\, \nat\big( \d_{\zb}\om(g,z) A_{\wb}(h,w) + A_{\zb}(g,z) \d_{\wb}\om(h,w) \big)  
$$
$$
= \int \frac{d^2z\, d^2w}{2\pi^2}\, \frac{1}{(hg(v)-v)}  \d_{\zb}\left( \frac{1}{hg(v)-z}+\frac{1}{z-v} \right) \frac{z-v}{g(z)-g(v)} \nat \om(g,z)A_{\wb}(h,w) $$
$$
+ \int\frac{ d^2z\, d^2w}{2\pi^2}\, \d_{\vb}\left( \frac{1}{(hg(v)-v)} \left( \frac{1}{hg(v)-z}+\frac{1}{z-v} \right) \right) \frac{z-v}{g(z)-g(v)} \nat A_{\zb}(g,z) \frac{\d\vb}{\d\wb} \om(h,w) 
$$
$$
= \int\frac{d^2z\, d^2w}{2\pi}\, \frac{1}{(hg(v)-v)} \left( -\delta^2(hg(v)-z)+\delta^2(z-v) \right) \frac{z-v}{g(z)-g(v)} \nat \om(g,z)A_{\wb}(h,w)  
$$
$$
+ \int\frac{d^2z\, d^2w}{2\pi}\, \frac{1}{(hg(v)-v)} \left( \frac{\d\overline{hg(v)}}{\d\vb} \delta^2(hg(v)-z) - \delta^2(z-v) \right) \frac{z-v}{g(z)-g(v)} \nat A_{\zb}(g,z) \frac{\d\vb}{\d\wb} \om(h,w) 
$$
$$
+ \int\frac{d^2z\, d^2w}{2\pi^2}\, \d_{\vb}\frac{1}{(hg(v)-v)} \left( \frac{1}{hg(v)-z}+\frac{1}{z-v} \right) \frac{z-v}{g(z)-g(v)} \nat A_{\zb}(g,z) \frac{\d\vb}{\d\wb} \om(h,w) \ .
$$
Before evaluating separately these last three integrals, we need to establish some identity. Around the fixed point $v_0$ of order $n$ (or equivalently the fixed point $w_0=g(v_0)$) one has the equality of distributions
$$
\d_v^{n-1}\left(\frac{1}{v-v_0}\right)\, \left(\frac{v-v_0}{w-w_0}\right)^n = \d_w^{n-1}\left(\frac{1}{w-w_0}\right) \ .
$$
The quotient $(hg(v)-v)/(gh(w)-w)$ is a holomorphic function of $v$ in the neighborhood of $v_0$ since $gh(w)-w=ghg(v)-g(v)$. Hence we deduce the \emph{distributional} identity
\beq
\frac{1}{hg(v)-v}\, \frac{hg(v)-v}{gh(w)-w} &=& \frac{(-)^{n-1}}{(n-1)!}\, \d_v^{n-1} \left(\frac{1}{v-v_0}\right) \, H^n_{hg,v_0}(v) \, \frac{hg(v)-v}{gh(w)-w}  \non\\
&=& \frac{(-)^{n-1}}{(n-1)!}\, \d_v^{n-1} \left(\frac{1}{v-v_0}\right) \, \frac{(v-v_0)^n}{hg(v)-v} \, \frac{hg(v)-v}{gh(w)-w} \non\\
&=& \frac{(-)^{n-1}}{(n-1)!}\, \d_v^{n-1} \left(\frac{1}{v-v_0}\right) \, \left(\frac{v-v_0}{w-w_0}\right)^n \, \frac{(w-w_0)^n}{gh(w)-w} \non\\
&=& \frac{(-)^{n-1}}{(n-1)!}\, \d_w^{n-1}\left(\frac{1}{w-w_0}\right) H^n_{gh,w_0}(w) \non\\
&=& \frac{1}{gh(w)-w} \ .\non
\eeq
Now we can evaluate the first integral:
\beq
\lefteqn{ \int \frac{d^2z\, d^2w}{2\pi}\, \frac{1}{(hg(v)-v)} \left( -\delta^2(hg(v)-z)+\delta^2(z-v) \right) \frac{z-v}{g(z)-g(v)} \nat \om(g,z)A_{\wb}(h,w) } \non\\
&=& \int\frac{d^2w}{2\pi}\, \frac{1}{(hg(v)-v)} \left( -\frac{hg(v)-v}{gh(w)-w}\, \nat \om(g,h(w)) A_{\wb}(h,w)  + \frac{\d v}{\d w}\, \nat \om(g,v)A_{\wb}(h,w) \right)  \non\\
&=& -\int\frac{d^2w}{2\pi}\, \frac{\nat \om(g,h(w)) A_{\wb}(h,w)}{gh(w)-w} + \int\frac{ d^2v}{2\pi}\, \frac{\nat \om(g,v)\d_{\vb}\overline{g(v)} A_{\overline{g}}(h,g(v))}{(hg(v)-v)}  \ .\non
\eeq
We proceed similarly with the second integral:
$$
\int\frac{d^2z\, d^2w}{2\pi}\, \frac{1}{(hg(v)-v)} \left( \frac{\d\overline{hg(v)}}{\d\vb} \delta^2(hg(v)-z) - \delta^2(z-v) \right) \frac{z-v}{g(z)-g(v)} \nat A_{\zb}(g,z) \frac{\d\vb}{\d\wb} \om(h,w) 
$$
$$
= \int\frac{d^2w}{2\pi}\, \frac{1}{(hg(v)-v)} \left( \frac{\d\overline{h(w)}}{\d\vb} \frac{hg(v)-v}{gh(w)-w} \, \nat A_{\overline{h}}(g,h(w)) \frac{\d\vb}{\d\wb} \om(h,w)  - \left| \frac{\d v}{\d w} \right|^2 \nat A_{\vb}(g,v) \om(h,w) \right) 
$$
$$
= \int \frac{d^2w}{2\pi}\, \frac{\nat \d_{\wb}\overline{h(w)} A_{\overline{h}}(g,h(w)) \om(h,w)}{gh(w)-w}   - \int\frac{d^2v}{2\pi}\, \frac{\nat A_{\vb}(g,v) \om(h,g(v))}{hg(v)-v} \ .
$$
For the third integral we write explicitly the renormalized form of $1/(hg(v)-v)$:
$$
\int \frac{d^2z\, d^2w}{2\pi^2}\, \d_{\vb}\frac{1}{(hg(v)-v)} \left( \frac{1}{hg(v)-z}+\frac{1}{z-v} \right) \frac{z-v}{g(z)-g(v)} \nat A_{\zb}(g,z) \frac{\d\vb}{\d\wb} \om(h,w) 
$$
$$
= \int\frac{d^2z\, d^2w}{2\pi^2}\, \d_{\vb}\left( \frac{(-)^{n-1}}{(n-1)!}\, \d_v^{n-1} \frac{1}{v-v_0} \, H^n_{hg,v_0}(v)  \right) \left( \frac{1}{hg(v)-z}+\frac{1}{z-v} \right)
$$
$$
\times\ \frac{z-v}{g(z)-g(v)} \nat A_{\zb}(g,z) \frac{\d\vb}{\d\wb} \om(h,w) 
$$
$$
= \int \frac{d^2z\, d^2w}{2\pi}\, \frac{(-)^{n-1}}{(n-1)!}\, \d_v^{n-1} \delta^2(v-v_0) \, H^n_{hg,v_0}(v) \left( \frac{1}{hg(v)-z}+\frac{1}{z-v} \right)
$$
$$
\times \  \frac{z-v}{g(z)-g(v)} \nat A_{\zb}(g,z) \frac{\d\vb}{\d\wb} \om(h,w) 
$$
It vanishes because we already know the identity
$$
\d_v^{n-1} \delta^2(v-v_0) \, \left( \frac{1}{hg(v)-z}+\frac{1}{z-v} \right) = 0 \ .
$$
Collecting these results we see that the contribution of an isolated fixed point to the variation $\dd W_R^2(A)$ is, at linear order in $A$,
\beq
\lefteqn{ -\frac{1}{2}\int d^2z\, d^2w\, \frac{\nat\big( \d_{\zb}\om(g,z) A_{\wb}(h,w) + A_{\zb}(g,z) \d_{\wb}\om(h,w)\big)}{\pi^2(h(w)-z)(g(z)-w)}  } \non\\
&&\quad = \int\frac{d^2w}{2\pi}\, \frac{1}{gh(w)-w}\, \nat\big( \om(h,w) \frac{\d\overline{h(w)}}{\d\wb} A_{\overline{h}}(g,h(w)) - A_{\wb}(h,w) \om(g,h(w)) \big) \non\\
&&\qquad + \int\frac{d^2v}{2\pi}\, \frac{1}{hg(v)-v} \, \nat\big( \om(g,v)\frac{\d\overline{g(v)}}{\d \vb} A_{\overline{g}}(h,g(v)) - A_{\vb}(g,v) \om(h,g(v)) \big) \ .\non
\eeq
Hence summing over $g,h\in G$ it is exactly canceled by the analogous term coming from $\dd W_R^1(A)$. We conclude there is no contribution of the isolated fixed points to $\Delta^1(\om,A)$. It remains to look at the non-isolated fixed points. Hence suppose that $hg(z)=z$ in a neighborhood of $z_0$. The renomalization performed in vi) yields
$$
\frac{1}{(h(w)-z)(g(z)-w)} = \d_z\frac{1}{z-v}\, \left(\frac{z-v}{g(z)-g(v)}\right) 
$$
with the change of variable $w=g(v)$. Since $hg(v)=v$ one has $v=h(w)$. Integrating by parts we find 
\beq
\lefteqn{-\frac{1}{2}\int d^2z\, d^2w\, \frac{\nat\big( \d_{\zb}\om(g,z) A_{\wb}(h,w) + A_{\zb}(g,z) \d_{\wb}\om(h,w) \big)}{\pi^2(h(w)-z)(g(z)-w)} } \non\\
&&= \int\frac{d^2z\, d^2w}{2\pi^2}\, \d_z\left(\d_{\zb}\frac{1}{z-v}\right) \, \frac{z-v}{g(z)-g(v)} \, \nat \om(g,z) A_{\wb}(h,w) \non\\
&&\quad + \int \frac{d^2z\, d^2w}{2\pi^2}\, \d_z\left(\d_{\vb}\frac{1}{z-v}\right) \, \frac{z-v}{g(z)-g(v)} \, \nat A_{\zb}(g,z) \frac{\d\vb}{\d\wb}\om(h,w)  \non\\
&&= -\int \frac{d^2z\, d^2w}{2\pi}\, \delta^2(z-v) \, \d_z\left( \frac{z-v}{g(z)-g(v)} \, \nat \om(g,z) A_{\wb}(h,w) \right) \non\\
&&\quad + \int \frac{d^2z\, d^2w}{2\pi}\, \delta^2(z-v) \, \d_z\left( \frac{z-v}{g(z)-g(v)} \, \nat A_{\zb}(g,z) \frac{\d\vb}{\d\wb}\om(h,w) \right) \non
\eeq
We Taylor expand $g(z)-g(v)$ around the diagonal $z-v=0$:
$$
\frac{z-v}{g(z)-g(v)} = \big( g'(v) + \frac{z-v}{2}g''(v) +\ldots \big)^{-1}\ .
$$
Then perform the integral over $z$ and recall that $g'(v)=\d w/\d v$ and $h'(w)=\d v/\d w$:
\beq
\lefteqn{ -\frac{1}{2}\int d^2z\, d^2w\, \frac{\nat\big( \d_{\zb}\om(g,z) A_{\wb}(h,w) + A_{\zb}(g,z) \d_{\wb}\om(h,w) \big)}{\pi^2(h(w)-z)(g(z)-w)}  } \non\\
&=& -\int \frac{d^2w}{2\pi}\, \left( -\frac{\frac{1}{2}g''(v)}{g'(v)^2} \, \nat \om(g,v) A_{\wb}(h,w) + \frac{1}{g'(v)} \, \nat \d_v\om(g,v) A_{\wb}(h,w) \right) \non\\
&& + \int\frac{d^2w}{2\pi}\, \left( -\frac{\frac{1}{2}g''(v)}{g'(v)^2}  \, \nat A_{\vb}(g,v) \frac{\d\vb}{\d_{\wb}}\om(h,w) + \frac{1}{g'(v)}\, \nat \d_v A_{\vb}(g,v) \frac{\d\vb}{\d_{\wb}}\om(h,w) \right) \non\\
&=& -\int\frac{d^2w}{2\pi}\, \frac{\d v}{\d w}\big( -\frac{1}{2}\d_v\ln g'(v) \, \nat \om(g,v) A_{\wb}(h,w) + \nat \d_v\om(g,v) A_{\wb}(h,w) \big) \non\\
&& + \int\frac{d^2w}{2\pi}\, \left|\frac{\d v}{\d w}\right|^2 \big( -\frac{1}{2}\d_v\ln g'(v) \, \nat A_{\vb}(g,v)\om(h,w) + \nat \d_v A_{\vb}(g,v) \om(h,w) \big) \non\\
&=& -\int\frac{d^2w}{2\pi}\, \big( \frac{1}{2}\d_w\ln h'(w) \, \nat\om(g,h(w)) A_{\wb}(h,w) + \nat \d_w\om(g,h(w)) A_{\wb}(h,w) \big) \non\\
&& + \int\frac{d^2v}{2\pi}\, \big( -\frac{1}{2}\d_v\ln g'(v) \, \nat A_{\vb}(g,v)\om(h,w) + \nat \d_v A_{\vb}(g,v) \om(h,w) \big) \non\\
&=& \int\frac{d^2w}{2\pi}\, \nat\big( \d_w -\frac{1}{2}\d_w\ln h'(w) \big) A_{\wb}(h,w) \, \om(g,h(w)) \non\\
&& +\int \frac{d^2v}{2\pi}\, \nat\big( \d_v -\frac{1}{2}\d_v\ln g'(v) \big) A_{\vb}(g,v) \, \om(h,g(v))\ .\non
\eeq
Hence summing over $g,h\in G$ we get the desired expression for the anomaly:
$$
\Delta^1(\om,A)=\frac{1}{\pi} \sum_{g,h\in G}\int_{z=hg(z)}d^2z\, \nat \big(\d_z- \frac{1}{2}\d_z\ln g'(z)\big)A_{\zb}(g,z)\,\,\om(h,g(z))\ .
$$ 
It remains to show that the degree two component of the anomaly $\Delta^2(\om,A)$ vanishes. We first rewrite the contribution of $\dd W_R^2(A)$ as
\beq
\lefteqn{ -\frac{1}{2}\sum_{g,h\in G}\int d^2z\, d^2w\, \frac{\nat\big( [A,\om]_{\zb}(g,z) A_{\wb}(h,w) + A_{\zb}(g,z) [A,\om]_{\wb}(h,w) \big) }{\pi^2(h(w)-z)(g(z)-w)} } \non\\
&=&\sum_{g,h,i\in G}\int \frac{d^2z\, d^2w}{-2\pi^2}\,  \frac{\nat \big( A_{\zb}(g,z)\om(i,g(z)) - \om(g,z) \frac{\d\overline{g(z)}}{\d\zb}A_{\gb}(i,g(z)) \big) A_{\wb}(h,w)}{(h(w)-z)(ig(z)-w)}  \non\\
&+& \sum_{g,h,i\in G}\int \frac{d^2z\, d^2w}{-2\pi^2}\,   \frac{ \nat A_{\zb}(g,z) \big( A_{\wb}(i,w)\om(h,i(w)) - \om(i,w) \frac{\d\overline{i(w)}}{\d\wb}A_{\ib}(h,i(w)) \big)}{(hi(w)-z)(g(z)-w)}  \non
\eeq
On the other hand, $\Delta^2(\om,A)$ gets also a contribution from the variation of $W^3(A)$:
$$
\dd W^3(A) = \sum_{g,h,i\in G} \int d^2z\, d^2w\, d^2v\, \frac{\nat \d_{\zb}\om(g,z) A_{\wb}(h,w) A_{\vb}(i,v)}{\pi^3(i(v)-z)(g(z)-w)(h(w)-v)}  + O(A^3)\ .
$$
Let us perform the change of variables $w=g(s)$, $v=hg(u)$ and rewrite the distribution kernel as
\beq
\lefteqn{ \frac{1}{(i(v)-z)(g(z)-w)(h(w)-v)} = } \non\\
&& \frac{z-s}{g(z)-g(s)} \, \left(\frac{1}{ihg(u)-z}+\frac{1}{z-s}\right) \frac{1}{(ihg(u)-s)(s-u)}\, \frac{s-u}{hg(s)-hg(u)} \ ,\non
\eeq
where $1/((ihg(u)-s)(s-u))$ is renormalized by iv) - vi). Then we integrate by parts, taking into account that $(z-s)/(g(z)-g(s))$ is a holomorphic function:
\beq
\lefteqn{ \int d^2z\, d^2w\, d^2v\, \frac{\nat \d_{\zb}\om(g,z) A_{\wb}(h,w) A_{\vb}(i,v)}{\pi^3(i(v)-z)(g(z)-w)(h(w)-v)} }\non\\
&=& \int \frac{d^2z\, d^2w\, d^2v}{\pi^2}\, \frac{z-s}{g(z)-g(s)}\, \left(\delta^2(ihg(u)-z)- \delta^2(z-s)\right) \non\\
&& \qquad \times  \frac{1}{(ihg(u)-s)(s-u)}\, \frac{s-u}{hg(s)-hg(u)}\, \nat \om(g,z) A_{\wb}(h,w) A_{\vb}(i,v) \non\\
&=& \int \frac{d^2w\, d^2v}{\pi^2}\, \frac{ihg(u)-s}{gihg(u)-g(s)}\,\frac{\nat \om(g,ihg(u)) A_{\wb}(h,w) A_{\vb}(i,v)}{(ihg(u)-s)(s-u)}\, \frac{s-u}{hg(s)-hg(u)} \non\\
&& - \int \frac{d^2w\, d^2v}{\pi^2}\, \left(\frac{\d g(s)}{\d s}\right)^{-1}\,\frac{\nat \om(g,s) A_{\wb}(h,w) A_{\vb}(i,v)}{(ihg(u)-s)(s-u)}\, \frac{s-u}{hg(s)-hg(u)} \non
\eeq
Now suppose $ihg$ has an isolated fixed point. The renormalization iv) - v) of the function $1/((ihg(u)-s)(s-u))$ implies
\beq
\lefteqn{\frac{ihg(u)-s}{gihg(u)-g(s)}\,\frac{1}{(ihg(u)-s)(s-u)} }\non\\
&&= \frac{1}{ihg(u)-u} \left( \frac{1}{gihg(u)-g(s)} + \frac{ihg(u)-s}{gihg(u)-g(s)}\, \frac{1}{s-u} \right) \ . \non
\eeq
The distribution in parenthesis is a locally integrable function of the variable $s$ (it has singularity order $-1$). On the other hand we have
\beq
\lefteqn{ \frac{1}{(gihg(u)-g(s))(s-u)} }\non\\
&=& \frac{1}{gihg(u)-g(u)} \left( \frac{1}{gihg(u)-g(s)} + \frac{1}{g(s)-g(u)} \right) \frac{g(s)-g(u)}{s-u} \non\\
&=& \frac{1}{ihg(u)-u}\frac{ihg(u)-u}{gihg(u)-g(u)}\left( \frac{1}{gihg(u)-g(s)} + \frac{1}{g(s)-g(u)} \right) \frac{g(s)-g(u)}{s-u}\ . \non
\eeq
But in the sense of locally integrable functions one has the equality
\beq
\lefteqn{\frac{ihg(u)-u}{gihg(u)-g(u)}\left( \frac{1}{gihg(u)-g(s)} + \frac{1}{g(s)-g(u)} \right) \frac{g(s)-g(u)}{s-u} } \non\\
&& \mspace{150mu} = \frac{1}{gihg(u)-g(s)} + \frac{ihg(u)-s}{gihg(u)-g(s)}\, \frac{1}{s-u}\ , \non
\eeq
which shows the (non-obvious) equality of distributions
$$
\frac{ihg(u)-s}{gihg(u)-g(s)}\,\frac{1}{(ihg(u)-s)(s-u)} = \frac{1}{(gihg(u)-g(s))(s-u)}\ ,
$$
and also
\beq
\lefteqn{\frac{ihg(u)-s}{gihg(u)-g(s)}\,\frac{1}{(ihg(u)-s)(s-u)}\,\frac{s-u}{hg(s)-hg(u)} } \label{bigre} \\
&& \mspace{150mu}  = \frac{1}{(gihg(u)-g(s))(hg(s)-hg(u))} \ . \non
\eeq
If $ihg$ is the identity around the fixed point, by renormalization vi) we have 
\beq
\frac{u-s}{g(u)-g(s)}\, \frac{1}{(u-s)(s-u)} &=& \frac{u-s}{g(u)-g(s)}\, \d_u\frac{1}{u-s} \non\\
&=& \d_{g(u)} \left(\frac{1}{g(u)-g(s)}\right) \, \frac{g(s)-g(u)}{s-u}\non\\
&=& \frac{1}{(g(u)-g(s))(g(s)-g(u))} \, \frac{g(s)-g(u)}{s-u}\non\\
&=& \frac{1}{(g(u)-g(s))(s-u)}\ , \non
\eeq
and (\ref{bigre}) is also valid. Now we can continue the computation
\beq
\lefteqn{ \int d^2z\, d^2w\, d^2v\, \frac{\nat \d_{\zb}\om(g,z) A_{\wb}(h,w) A_{\vb}(i,v)}{\pi^3(i(v)-z)(g(z)-w)(h(w)-v)} }\non\\
&=& \int \frac{d^2w\, d^2v}{\pi^2}\, \frac{1}{(gihg(u)-g(s))(hg(s)-hg(u))} \, \nat \om(g,ihg(u)) A_{\wb}(h,w) A_{\vb}(i,v) \non\\
&& - \int \frac{d^2s\, d^2v}{\pi^2}\, \frac{\d \overline{g(s)}}{\d \sb}\, \frac{1}{(ihg(u)-s)(hg(s)-hg(u))} \, \nat \om(g,s) A_{\wb}(h,w) A_{\vb}(i,v) \non\\
&=& \int \frac{d^2w\, d^2v}{\pi^2}\, \frac{1}{(gi(v)-w)(h(w)-v)} \, \nat \om(g,i(v)) A_{\wb}(h,w) A_{\vb}(i,v) \non\\
&& - \int \frac{d^2s\, d^2v}{\pi^2}\, \frac{1}{(i(v)-s)(hg(s)-v)} \, \nat \om(g,s) \frac{\d\overline{g(s)}}{\d \sb}  A_{\gb}(h,g(s)) A_{\vb}(i,v) \non
\eeq
Then summing over $g,h,i\in G$ and renaming the variables we obtain
\beq
\lefteqn{\frac{1}{2}\sum_{g,h,i} \int d^2z\, d^2w\, d^2v\, \frac{\nat \d_{\zb}\om(g,z) A_{\wb}(h,w) A_{\vb}(i,v)}{\pi^3(i(v)-z)(g(z)-w)(h(w)-v)} }\non\\
&=& \sum_{g,h,i} \int \frac{d^2z\, d^2w}{2\pi^2}\, \frac{1}{(hi(w)-z)(g(z)-w)} \, \om(h,i(w)) A_{\zb}(g,z) A_{\wb}(i,w) \non\\
&& - \sum_{g,h,i} \int \frac{d^2w\, d^2z}{2\pi^2}\, \frac{1}{(h(w)-z)(ig(z)-w)} \, \nat \om(g,z) \frac{\d\overline{g(z)}}{\d \zb}  A_{\gb}(i,g(z)) A_{\wb}(h,w) \non
\eeq
This quantity exactly compensates half of the terms in the contribution of $\dd W_R^2(A)$. For the other half, we rewrite the distribution kernel with the same change of variables $w=g(s)$, $v=hg(u)$:
\beq
\lefteqn{ \frac{1}{(i(v)-z)(g(z)-w)(h(w)-v)} = } \non\\
&& \frac{s-u}{hg(s)-hg(u)} \, \frac{1}{(s-u)(ihg(u)-s)}\, \left(\frac{1}{ihg(u)-z}+\frac{1}{z-s}\right) \frac{z-s}{g(z)-g(s)} \ ,\non
\eeq
where $1/((s-u)(ihg(u)-s))$ is renormalized by iv) - vi). Then integrate by parts:
\beq
\lefteqn{ \int d^2z\, d^2w\, d^2v\, \frac{\nat \d_{\zb}\om(g,z) A_{\wb}(h,w) A_{\vb}(i,v)}{\pi^3(i(v)-z)(g(z)-w)(h(w)-v)} }\non\\
&=&  \int \frac{d^2z\, d^2w\, d^2v}{\pi^2}\, \frac{s-u}{hg(s)-hg(u)} \, \frac{1}{(s-u)(ihg(u)-s)} \non\\
&& \times \left(\delta^2(ihg(u)-z) - \delta^2(z-s)\right) \frac{z-s}{g(z)-g(s)} \, \nat \om(g,z) A_{\wb}(h,w) A_{\vb}(i,v) \non\\
&=&  \int \frac{d^2w\, d^2v}{\pi^2}\, \frac{s-u}{hg(s)-hg(u)} \, \frac{\nat \om(g,ihg(u)) A_{\wb}(h,w) A_{\vb}(i,v)}{(s-u)(ihg(u)-s)}\, \frac{ihg(u)-s}{gihg(u)-g(s)} \non\\
&& - \int \frac{d^2w\, d^2v}{\pi^2}\, \frac{s-u}{hg(s)-hg(u)} \, \frac{\nat \om(g,s) A_{\wb}(h,w) A_{\vb}(i,v)}{(s-u)(ihg(u)-s)}\, \left(\frac{\d g(s)}{\d s}\right)^{-1}\ . \non
\eeq
Use the same tricks leading to (\ref{bigre}) and rewrite this as
\beq
\lefteqn{ \int d^2z\, d^2w\, d^2v\, \frac{\nat \d_{\zb}\om(g,z) A_{\wb}(h,w) A_{\vb}(i,v)}{\pi^3(i(v)-z)(g(z)-w)(h(w)-v)} }\non\\
&=&  \int \frac{d^2w\, d^2v}{\pi^2}\, \frac{1}{(hg(s)-hg(u))(gihg(u)-g(s))} \, \nat \om(g,ihg(u)) A_{\wb}(h,w) A_{\vb}(i,v) \non\\
&& - \int \frac{d^2s\, d^2v}{\pi^2}\, \frac{1}{(hg(s)-hg(u))(ihg(u)-s)}\, \frac{\d \overline{g(s)}}{\d \sb}\, \nat \om(g,s) A_{\wb}(h,w) A_{\vb}(i,v) \non\\
&=&  \int \frac{d^2w\, d^2v}{\pi^2}\, \frac{1}{(h(w)-v)(gi(v)-w)} \, \nat \om(g,i(v)) A_{\wb}(h,w) A_{\vb}(i,v) \non\\
&& - \int \frac{d^2s\, d^2v}{\pi^2}\, \frac{1}{(hg(s)-v)(i(v)-s)}\, \nat \om(g,s) \frac{\d \overline{g(s)}}{\d \sb} A_{\gb}(h,g(s)) A_{\vb}(i,v) \non
\eeq
Summing over $g,h,i\in G$ and renaming the variables we obtain
\beq
\lefteqn{\frac{1}{2}\sum_{g,h,i} \int d^2z\, d^2w\, d^2v\, \frac{\nat \d_{\zb}\om(g,z) A_{\wb}(h,w) A_{\vb}(i,v)}{\pi^3(i(v)-z)(g(z)-w)(h(w)-v)} }\non\\
&=& \sum_{g,h,i} \int \frac{d^2w\, d^2z}{2\pi^2}\, \frac{1}{(h(w)-z)(ig(z)-w)} \, \nat \om(i,g(z)) A_{\wb}(h,w) A_{\zb}(g,z) \non\\
&& - \sum_{g,h,i} \int \frac{d^2w\, d^2z}{2\pi^2}\, \frac{1}{(hi(w)-z)(g(z)-w)}\, \nat \om(i,w) \frac{\d \overline{i(w)}}{\d \wb} A_{\ib}(h,i(w)) A_{\zb}(g,z) \non
\eeq
which compensates the other half of terms in the contribution of $\dd W_R^2(A)$. Hence $\Delta^2(\om,A)=0$. \cqfd \\

The contribution of an isolated fixed point of order $n$ in the degree zero part of the anomaly $\Delta^0(\om,A)$ involves only the derivatives of the mapping $g$ up to order $2n-1$. For example in low orders one finds
\beq
\lefteqn{\frac{-1}{(n-1)!}\, \d_z^{n-1}\big(H_{g,z_0}^n(z) \nat\om(g,z)\big)_{z=z_0} = } \non\\
&(n=1)&:\qquad \frac{1}{1-g'(z_0)}\, \nat\om(z_0)    \non\\
&(n=2)&:\qquad \frac{2}{g''(z_0)} \left( \frac{1}{3} \frac{g'''(z_0)}{g''(z_0)} \, \nat\om(z_0) -  \nat\d_z\om(z_0) \right)   \non\\
&(n=3)&:\qquad \frac{3}{2g'''(z_0)} \left( \frac{1}{10} \frac{g^{(5)}(z_0)}{g'''(z_0)}\right.  \nat\om(z_0) - \frac{1}{8} \left(\frac{g^{(4)}(z_0)}{g'''(z_0)}\right)^2 \, \nat\om(z_0)  \non\\
&& \qquad \qquad \qquad \qquad \left. + \frac{1}{2} \frac{g^{(4)}(z_0)}{g'''(z_0)} \, \nat \d_z\om(z_0) - \nat \d_z^2\om(z_0) \right)   \non
\eeq
One recovers the well-known Lefschetz numbers in the case $n=1$, whereas for $n>1$ the higher order jets of the mapping $g$ are involved.

\section{Index theorem}\label{sind}

In this section we shall use the anomaly formula established in Proposition \ref{ploc} to calculate the diagonal of the commutative diagram (\ref{rrg}). Thus as before we let $\Si$ be the complex plane and $G$ be a discrete group acting by conformal transformations on $\Si$. Consider the set
\be
\Gamma = \coprod_{g\in G} \dom(g) = \{(g,z)\in G\times \Si\ | \ z\in\dom(g)\}\ .
\ee
It is a smooth \'etale groupoid for the composition law $(g,z)\cdot(h,w) = (gh,w)$ whenever $z=h(w)$, its set of units corresponding to $\Gamma_0=\dom(1)\subset\Si$. The crossed product $\Ac_0=\cinfc(\Si)\cp G$ co\"{\i}ncides with the convolution algebra of smooth, compactly supported functions over $\Gamma$: to an element $a=\sum_{g\in G}a(g) U^*_g$ of $\Ac_0$ corresponds the function $(g,z)\mapsto a(g,z)=a(g)(z)$ over the groupoid. \\
We say that a morphism $\gamma_0=(g,z_0)\in \Gamma$ is an automorphism if $z_0=g(z_0)$ is a fixed point. The order of $\gamma_0$ is the integer $n\in \nn\cup\{\infty\}$ corresponding to the order of the fixed point. According to the discussion of section \ref{sren}, an automorphism of order $n<\infty$ is necessarily isolated, whereas for $n=\infty$ all morphisms are automorphisms in the vicinity of $\gamma_0$. Hence the set of automorphisms is the union of a discrete set $\Gamma_f$ (for $n<\infty$) and a one-dimensional complex manifold $\Gamma_{\infty}$ (for $n=\infty$). $\Gamma_{\infty}$ contains the set of units $\Gamma_0$. Proposition \ref{pano} precisely shows that the anomaly splits as the sum of $\Delta^0(\om,A)$ localized at $\Gamma_f$, and $\Delta^1(\om,A)$ localized at $\Gamma_{\infty}$. Let us describe these terms in a more intrinsic way.\\

The component of degree zero in the anomaly (\ref{d0}) is a sum over all isolated automorphisms $\gamma_0=(g,z_0)$ of the numbers
$$
\frac{-1}{(n-1)!}\, \d_z^{n-1} \big( H^n_{g,z_0}(z)\,\nat\om(g,z) \big)_{z=z_0}\ ,\qquad H^n_{g,z_0}(z)=\frac{(z-z_0)^n}{g(z)-z}\ ,
$$
where $n<\infty$ is the order of $\gamma_0$. Remark that replacing $n$ by any integer $m\geq n$ does not affect the result. The above quantity a priori depends on the complex coordinate system $z$, and one may wonder if it can be defined intrinsically, i.e. by means of the complex structure on $\Si$ only.

\begin{lemma}\label{linv}
Let $\gamma_0=(g,z_0)\in\Gamma_f$ be an isolated automorphism of order $n$. The linear functional $\Ac_0\to\cc$ given by 
$$
a\mapsto \d_z^{n-1} \big( H^n_{g,z_0}(z)\,a(g,z) \big)_{z=z_0}
$$
is intrinsically defined, i.e. independent of the complex coordinate system $z$. We write $\gamma=(g,z)\in \Gamma$ in a neighborhood of $\gamma_0$ endowed with its complex structure, $H^n_{g,z_0}(z)=H^n_{\gamma_0}(\gamma)$ and identify $\d_z$ with $\d_{\gamma}$. Then summing over all isolated automorphisms the linear functional $\Phi(\Gamma):\Ac_0 \to \cc$
$$
\Phi(\Gamma)(a)= \sum_{\gamma_0 \in\Gamma_f} \frac{-1}{(n-1)!}\, \d_{\gamma}^{n-1} \big( H^n_{\gamma_0}(\gamma)\,a(\gamma) \big)_{\gamma=\gamma_0}
$$
is a trace on $\Ac_0$.
\end{lemma}
{\it Proof:} Let $(g,z_0)$ be an isolated automorphism of order $n$, and $f\in \cinfc(\Si)$ be a test function. We compute
\beq
\d_z^{n-1} \big( H^n_{g,z_0}(z) f(z) \big)_{z=z_0} &=& (-)^{n-1}\int d^2z\, \d_z^{n-1}\delta^2(z-z_0) \frac{(z-z_0)^n}{g(z)-z} f(z) \non\\
&=& (-)^{n}\int d^2z\, \d_z^{n-1}\frac{1}{\pi(z-z_0)} \frac{(z-z_0)^n}{g(z)-z} \d_{\zb}f(z) \ . \non
\eeq
Let $w=h(z)$ be a change of complex coordinate. Hence $h$ is a holomorphic function with $h'=\d w/\d z \neq 0$, and the point $w_0=h(z_0)$ is fixed by the conformal mapping $hgh^{-1}$. We know the distributional identity
$$
\d_z^{n-1} \frac{1}{z-z_0} = \frac{(w-w_0)^n}{(z-z_0)^n}\, \d_w^{n-1} \frac{1}{w-w_0}\ ,
$$
so making a change of variables in the integral yields
\beq
\lefteqn{\d_z^{n-1} \big( H^n_{g,z_0}(z) f(z) \big)_{z=z_0} = (-)^{n}\int d^2z\, \d_w^{n-1}\frac{1}{\pi(w-w_0)} \frac{(w-w_0)^n}{g(z)-z} \d_{\zb}f(z) } \non\\
&& \qquad =  (-)^{n}\int d^2w\, \frac{\d z}{\d w} \d_w^{n-1}\frac{1}{\pi(w-w_0)} \frac{(w-w_0)^n}{g(z)-z} \d_{\wb}f(h^{-1}(w)) \non\\
&& \qquad =  (-)^{n-1}\int d^2w\, \frac{\d z}{\d w} \d_w^{n-1}\delta^2(w-w_0) \frac{(w-w_0)^n}{g(z)-z} f(h^{-1}(w)) \non\\
&& \qquad = \d_w^{n-1} \Big( \frac{\d z}{\d w} \frac{(w-w_0)^n}{g(z)-z} f(h^{-1}(w)) \Big)_{w=w_0} \ . \non
\eeq
Now remark that $hgh^{-1}(w)-w = hg(z)-h(z)$ so we can write
$$
\frac{\d z}{\d w} \frac{(w-w_0)^n}{g(z)-z} = \frac{\d z}{\d w} \frac{hg(z)-h(z)}{g(z)-z} \frac{(w-w_0)^n}{hgh^{-1}(w)-w} \ .
$$
The difference $g(z)-z$ is of order $(z-z_0)^n$, hence
$$
\frac{\d z}{\d w} \frac{hg(z)-h(z)}{g(z)-z} = \frac{\d z}{\d w} \Big( \frac{\d h(z)}{\d z} + O((z-z_0)^n) \Big) = 1 + O((z-z_0)^n) \ .
$$
Equivalently it is $1 +O((w-w_0)^n) $. Differentiating $n-1$ times with respect to $w$ and taking the value at $w_0$ eliminates the remainder so one gets
\beq
\d_z^{n-1} \big( H^n_{g,z_0}(z) f(z) \big)_{z=z_0} &=& \d_w^{n-1} \Big( \frac{(w-w_0)^n}{hgh^{-1}(w)-w} f(h^{-1}(w)) \Big)_{w=w_0} \non\\
&=& \d_w^{n-1} \big( H^n_{hgh^{-1},w_0}(w) f(h^{-1}(w)) \big)_{w=w_0} \ . \non
\eeq
Hence this quantity does not depend on the choice of complex coordinate system. Now take $a,b\in \Ac_0$ and compute
\beq
\Phi(ab) &=& \sum_{\substack{g\in G \\ z_0=g(z_0)}} \frac{-1}{(n-1)!}\d_z^{n-1}\big( H^n_{g,z_0}(z) ab(g,z) \big)_{z=z_0} \non\\
&=& \sum_{\substack{g,h\in G \\ z_0=g(z_0)}} \frac{-1}{(n-1)!}\d_z^{n-1}\big( H^n_{g,z_0}(z) a(h,z)b(gh^{-1},h(z)) \big)_{z=z_0} \ . \non
\eeq
Change $g$ to $gh$ and rewrite this as
$$
\Phi(ab)= \sum_{\substack{g,h\in G \\ z_0=gh(z_0)}} \frac{-1}{(n-1)!}\d_z^{n-1}\big( H^n_{gh,z_0}(z) a(h,z)b(g,h(z)) \big)_{z=z_0} \ .
$$
Make the change of variables $w=h(z)$, so $w_0=h(z_0)$ is fixed by $hg$. The formula established above gives
\beq
\lefteqn{ \d_z^{n-1}\big( H^n_{gh,z_0}(z) a(h,z)b(g,h(z)) \big)_{z=z_0} }\non\\
&& \mspace{100mu} = \d_w^{n-1}\big( H^n_{hg,w_0}(w) a(h,h^{-1}(w))b(g,w) \big)_{w=w_0} \non \\
&& \mspace{100mu} = \d_w^{n-1}\big( H^n_{hg,w_0}(w) b(g,w) a(h,g(hg)^{-1}(w)) \big)_{w=w_0} \ .\non 
\eeq
By hypothesis $w_0$ is a fixed point of order $n$ for $hg$, hence $g(hg)^{-1}(w)$ differs from $g(w)$ by a term of order $(w-w_0)^n$ which will disappear after taking $n-1$ derivatives at the value $w=w_0$. Thus one has
$$
\Phi(ab) = \sum_{\substack{g,h\in G \\ w_0=hg(w_0)}} \frac{-1}{(n-1)!} \d_w^{n-1}\big( H^n_{hg,w_0}(w) b(g,w) a(h,g(w)) \big)_{w=w_0} = \Phi(ba)\ , 
$$
showing that $\Phi$ is a trace on $\Ac_0$. \cqfd\\

Recall that we decomposed the Maurer-Cartan form $\om_{\pm} = \sum_{g\in G} \om(g) r(g)_{\pm}$ by means of the components $\om(g)\in M_{\infty}(\cinfc(\Si))\otimes\Om^1\Th\Bc$. Replacing the representations $r(g)_{\pm}$ by the universal symbol $U^*_g$, we may view the Maurer-Cartan form as an element
\be
\omt=\sum_{g\in G}\om(g)U^*_g \ \in M_{\infty}(\Ac_0)\otimes\Om^1\Th\Bc\ ,
\ee
where the space $M_{\infty}(\Ac_0)\otimes\Om^1\Th\Bc$ is naturally a bimodule over the pro-algebra $M_{\infty}(\Ac_0)\otimes \Th\Bc$. Thus tensoring $\Phi(\Gamma):M_{\infty}(\Ac_0)\to\cc$ with the quotient map $\nat:\Om^1\Th\Bc\to \Om^1\Th\Bc_{\nat}$ yields a trace $\Phi(\Gamma)\nat$ on the bimodule $M_{\infty}(\Ac_0)\otimes\Om^1\Th\Bc$. The degree zero component of the anomaly is therefore
\be
\Delta^0(\om,A)= \Phi(\Gamma)\nat(\omt)\ .
\ee

The degree one component of the anomaly $\Delta^1(\om,A)$ is localized at the manifold of non-isolated automorphisms, so we have to introduce the convolution algebra of differential forms over the groupoid $\Gamma$. Let $\Omc^*(\Si)$ be the (bi)graded algebra of compactly supported differential forms over $\Si$. Since $G$ acts on $\Si$ by conformal mappings, it acts on $\Omc^*(\Si)$ by pullback, and the crossed product $\Omc^*(\Si)\cp G$ defines a (bi)graded algebra. $\Ac_0$ corresponds to its degree zero subalgebra. Moreover the differentials $\d=dz\d_z$ and $\deb=d\zb\d_{\zb}$ on $\Omc^*(\Si)$ commute with conformal transformations hence extend to differentials of bidegree $(1,0)$ and $(0,1)$ respectively on $\Omc^*(\Si)\cp G$. Of course $\d\deb+\deb\d=0$, and the de Rham operator $d=\d+\deb$ is another differential. Finally, there is a fourth differential coming from the action of $G$ on $\Si$: define the \emph{modular derivative} on $\Omc^*(\Si)\cp G$ as
\be
D(fU^*_g) =  \ln |g'|^2 f U^*_g\ ,\qquad \forall\ f\in \Omc^*(\Si)\ ,\ g\in G\ ,
\ee
where the scalar function $z\mapsto |g'(z)|^2$ measures the volume dilatation induced by the conformal mapping $g$. The chain rule immediately implies that $D$ is a derivation of degree zero on the algebra $\Omc^*(\Si)\cp G$. Hence by taking the commutator with $\d$ one obtains a differential 
\be
\delta=[\d,D]\ ,\qquad \delta^2=0\ ,
\ee
which anticommutes with $\d$, $\deb$ and $d$. Explicitly it reads $\delta(fU^*_g)= (\d \ln g')f U^*_g$. The manifold $\Gamma_{\infty}$ has a complex structure, therefore may be considered as an oriented two-dimensional real manifold. It thus defines a fundamental class on which genuine two-forms can be integrated. 

\begin{lemma}
The trilinear functional $[\Gamma] : \Ac_0\times\Ac_0\times\Ac_0\to\cc$ defined by
$$
[\Gamma] (a_0,a_1,a_2) = \int_{\Gamma_{\infty}} a_0da_1da_2
$$
is a cyclic 2-cocycle over $\Ac_0$ called the fundamental class of the groupoid $\Gamma$. The trilinear functional $c_1(\Gamma) : \Ac_0\times\Ac_0\times\Ac_0\to\cc$ defined by
$$
c_1(\Gamma) (a_0,a_1,a_2) = \int_{\Gamma_{\infty}} a_0(da_1\delta a_2 + \delta a_1 da_2)
$$
is a cyclic 2-cocycle over $\Ac_0$ called the first Chern class of the groupoid $\Gamma$. We define the Todd class as the cyclic 2-cocycle 
\be
\Td(\Gamma) :=  [\Gamma] -\frac{1}{2} c_1(\Gamma)\ .
\ee
\end{lemma}
{\it Proof:} See \cite{P0}. The fundamental class and the first Chern class are cyclic cocycles because integration over the manifold $\Gamma_{\infty}$ is a graded trace over the algebra $\Omc^*(\Si)\cp G$, and is closed for the differentials $d$ and $\delta$ (the modular derivative vanishes when localized at automorphisms). \cqfd\\

The Todd class has a simple expression using the differential $\nabla=d-\frac{1}{2}\delta$, $\nabla^2=0$. Since $\delta a_1\delta a_2=0$ for dimensional reasons, one gets
\be
\Td(\Gamma) (a_0,a_1,a_2)=\int_{\Gamma_{\infty}} a_0\nabla a_1 \nabla a_2\ .
\ee
In fact this expression is not a conformal invariant on $\Si$ because the modular derivative $D$ measures the dilatation $|g'|^2$ of the \emph{euclidean} volume $d\zb\wedge dz/2i$ by a mapping $g$, and thus depends on the choice of a volume form besides the complex structure. There is a priori no reason to prefer the euclidean volume and we could as well choose any smooth volume form $\nu$. The complex structure plus the volume form is equivalent to fixing a K\"ahler metric on $\Si$. As shown in \cite{P0}, the new modular derivative $D^{\nu}$ associated to $\nu$ modifies the differential $\delta^{\nu}=[\d,D^{\nu}]$, so that the square of $\nabla^{\nu}=d-\frac{1}{2} \delta^{\nu} $ does no longer vanish but is proportional to the curvature of the K\"ahler metric. Consequently the cyclic cocycle representing the first Chern class must be modified by adding a term proportional to the curvature, while its cyclic cohomology class remains unchanged. This establishes the link with the Todd class of ordinary (commutative) Riemann surfaces.\\

Turning to the gauge potential, we decomposed $A=\sum_{g\in G} A(g) r(g)_+$ by means of the one-forms $A(g)=d\zb A_{\zb}(g) \in M_{\infty}(\Omc^{0,1}(M))\otimes \Th\Bc$. Hence replacing $r(g)_+$ with the universal symbol $U^*_g$, we may view the potential as an element 
\be
\At=\sum_{g\in G} A(g)U^*_g \ \in M_{\infty}(\Omc^*(\Si)\cp G)\otimes \Th\Bc
\ee
Since $\Ac_0\subset \Omc^*(\Si)\cp G$ we interpret also $\omt$ as an element of $M_{\infty}(\Omc^*(\Si)\cp G)\otimes \Om^1\Th\Bc$, viewed is as bimodule over the DG pro-algebra $M_{\infty}(\Omc^*(\Si)\cp G)\otimes \Th\Bc$. The degree one component of the anomaly (\ref{d1}) then reads
\be
\Delta^1(\om,A)= \frac{-1}{2\pi i} \int_{\Gamma_{\infty}} \nat \nabla \At\, \omt\ .
\ee
These very simple expressions for $\Delta(\om,A)$ allow to calculate explicitly the diagonal of diagram (\ref{rrg})
$$
\vcenter{\xymatrix{
\Kt_i(\Ic\hotimes\Ac) \ar[r]^{\rho_!} \ar[d] \ar@{.>}[rd] & \Kt_i(\Ic\hotimes\Bc) \ar[d] \\
HP_i(\Ac) \ar[r]_{\ch(\rho)} & HP_i(\Bc) }}  \qquad i\in \zz_2\ ,
$$
acting on a class $[u]\in \Kt_1(\Ic\hotimes\Ac)$ of odd degree represented by an invertible $u\in M_{\infty}(\Ac_0)^+$, or a class $[e]\in \Kt_0(\Ic\hotimes\Ac)$ of even degree represented by an idempotent $e\in M_{\infty}(\Ac_0)$. First define a homomorphism $\tilde{\rho}: \Ac_0\to \Ac_0\otimes\Bc$ by $\tilde{\rho}(fU^*_g) = fU^*_g\otimes U^*_g$. Its lifting is a homomorphism $\tilde{\rho}_*: T\Ac_0 \to \Ac_0\otimes T\Bc$:
$$
\tilde{\rho}_*(f_1U^*_{g_1} \otimes\ldots\otimes f_kU^*_{g_k}) = \big( f_1f_2^{g_1} \ldots f_k^{g_{k-1}\ldots g_1} U^*_{g_k\ldots g_1} \big) \otimes (U^*_{g_1}\otimes \ldots \otimes U^*_{g_k})\ .
$$
$\tilde{\rho}_*$ extends to a homomorphism of pro-algebras $\Th\Ac_0 \to \Ac_0\otimes\Th\Bc$. In the odd case, take the canonical invertible lifting $\uh\in M_{\infty}(\Th\Ac_0)^+$ of $u$ and set 
\be
\ut=\tilde{\rho}_* \uh \ \in M_{\infty}(\Ac_0) \otimes \Th\Bc\ .
\ee 
Then taking $\omt = \ut^{-1}\dd\ut \in M_{\infty}(\Ac_0)\otimes\Om^1\Th\Bc$ as Maurer-Cartan form and $\At=\ut^{-1}\deb\ut \in M_{\infty}(\Omc^*(M)\cp G)\otimes\Th\Bc$ as gauge potential, the Chern character of $\rho_!(u)$ in periodic cyclic homology $HP_1(\Bc)$ is represented by the anomaly $\Delta(\om,A)$ up to a factor $\sqrt{2\pi i}$ (Proposition \ref{pano}). The result can be formulated nicely in terms of the Chern character of $\ut$ in the non-commutative de Rham homology of the algebra $\Ac_0\otimes \Th\Bc$:
$$
\ch_{\textup{dR}}(\ut) = \sum_{n\geq 0} \frac{n!}{(2n+1)!} \, \Tr\nat \left( \frac{\ut^{-1}\dd\ut}{\sqrt{2\pi i}} \right)^{2n+1}\ \in \Omh (\Ac_0\otimes \Th\Bc)_{\nat}\ .
$$
This is a cycle of odd degree in the complex $\Omh (\Ac_0\otimes \Th\Bc)_{\nat}=\Omh (\Ac_0\otimes \Th\Bc)/[\ ,\ ]$ gifted with the boundary $\dd$. In the even case, take the canonical idempotent lifting $\eh\in M_{\infty}(\Th\Ac_0)$ of $e$ and set 
\be
\et= \tilde{\rho}_*\eh \ \in M_{\infty}(\Ac_0)\otimes \Th\Bc\ . 
\ee
The Chern character of $\rho_!(e)$ in periodic cyclic homology $HP_0(\Bc)$ can also be reduced to an anomaly, using Bott periodicity. Again this can be expressed via the Chern character of $\et$ in non-commutative de Rham homology:
$$
\ch_{\textup{dR}}(\et) = \sum_{n\geq 0} \frac{(-)^n}{n!} \, \Tr\nat \left( \frac{\et\dd\et\dd\et}{2\pi i} \right)^n\ \in \Omh (\Ac_0\otimes \Th\Bc)_{\nat}\ .
$$
In fact $\ch_{\textup{dR}}(\ut)$ and $\ch_{\textup{dR}}(\et)$ represent the Chern characters of $\tilde{\rho}(u)$ and $\tilde{\rho}(e)$ in the cyclic homology of the discrete algebra $\Ac_0\otimes\Bc$. Now if $\varphi$ is any cyclic cocycle over $\Ac_0$, it induces a cap-product 
\be
\varphi \cap : HP_*(\Ac_0\otimes\Bc) \to HP_*(\Bc)
\ee
by composing the natural morphism of differential complexes $\Omh(\Ac_0\otimes \Th\Bc)_{\nat} \to (\Omh\Ac_0)_{\nat} \otimes (\Omh\Th\Bc)_{\nat}$ with the cyclic cocycle $\varphi: \Omh\Ac_0\to\cc$. The anomaly formula then shows that the diagonal map is nothing else but the cap-product of the Chern characters $\ch(\tilde{\rho}(u))\in HP_1(\Ac_0\otimes\Bc)$ and $\ch(\tilde{\rho}(e))\in HP_0(\Ac_0\otimes\Bc)$ with the cyclic cocycle $\varphi = \Phi(\Gamma) + \Td(\Gamma)$:

\begin{theorem}\label{tcup}
For any class $[u]\in \Kt_1(\Ic\hotimes\Ac)$ represented by an invertible $u\in M_{\infty}(\Ac_0)^+$, the Chern character of its pushforward $\rho_!(u)\in \Kt_1(\Ic\hotimes\Bc)$ is the periodic cyclic homology class of odd degree given by a cap-product
\beq
\ch(\rho_!(u)) &=& (\Phi(\Gamma)+ \Td(\Gamma))\cap \ch(\tilde{\rho}(u))  \\
&=& \Phi(\Gamma) \nat \Big(\frac{\ut^{-1}\dd\ut}{\sqrt{2\pi i}}\Big) -  \int_{\Gamma_{\infty}} \nat\, \frac{\ut^{-1}\nabla\ut\nabla\ut^{-1}\dd\ut}{2(2\pi i)^{3/2}} \ \in \Om^1\Th\Bc_{\nat}\ . \non
\eeq
For any class $[e]\in \Kt_0(\Ic\hotimes\Ac)$ represented by an idempotent $e\in M_{\infty}(\Ac_0)$, the Chern character of its pushforward $\rho_!(e)\in \Kt_0(\Ic\hotimes\Bc)$ is the periodic cyclic homology class of even degree given by a cap-product
\beq
\ch(\rho_!(e)) &=& (\Phi(\Gamma)+ \Td(\Gamma))\cap \ch(\tilde{\rho}(e)) \\
&=& \Phi(\Gamma) \nat (\et) -  \int_{\Gamma_{\infty}} \nat\, \frac{\et\nabla\et\nabla\et}{2\pi i}  \ \in \Th\Bc_{\nat}\ . \non
\eeq
\end{theorem}
{\it Proof:} For notational convenience suppose $u\in \Ac_0^+$. One has $\omt = \ut^{-1}\dd\ut$ and $\At=\ut^{-1}\deb\ut$. Let us rewrite the anomaly $\nabla \At\, \omt= \nabla(\ut^{-1}\deb\ut)\ut^{-1}\dd\ut$ in two different ways:
\beq
\nabla(\ut^{-1}\deb\ut)\ut^{-1}\dd\ut &=& \nabla(\ut^{-1}\deb\ut\ut^{-1})\dd\ut + \ut^{-1}\deb\ut\nabla\ut^{-1}\dd\ut \non \\
&=& -\nabla\deb \ut^{-1} \dd\ut + \ut^{-1}\deb\ut\nabla\ut^{-1}\dd\ut \non
\eeq
\beq
\nabla(\ut^{-1}\deb\ut)\ut^{-1}\dd\ut &=& \nabla\ut^{-1}\deb\ut\ut^{-1}\dd\ut + \ut^{-1}\nabla\deb\ut\ut^{-1}\dd\ut \non \\
&=& \ut^{-1}\nabla\ut\deb\ut^{-1}\dd\ut + \ut^{-1}\nabla\deb\ut\ut^{-1}\dd\ut \ .\non
\eeq
We sum these two equalities and apply the trace $\int_{\Gamma_{\infty}} \nat: (\Omc^*(M)\cp G)\otimes\Th\Bc \to \Th\Bc_{\nat}$, which is closed for the differentials $\nabla$ and $\deb$:
\beq
\lefteqn{2\int_{\Gamma_{\infty}}\nat\nabla(\ut^{-1}\deb\ut)\ut^{-1}\dd\ut = \int_{\Gamma_{\infty}}\nat(-\nabla\deb \ut^{-1} \dd\ut + \nabla\deb\ut\ut^{-1}\dd\ut\ut^{-1}) } \non\\
&& \mspace{120mu}  + \int_{\Gamma_{\infty}}\nat(\ut^{-1}\deb\ut\nabla\ut^{-1}\dd\ut + \ut^{-1}\nabla\ut\deb\ut^{-1}\dd\ut ) \non\\
&=& \int_{\Gamma_{\infty}}\nat(- \ut^{-1} \dd\nabla\deb \ut - \nabla\deb\ut\dd\ut^{-1})  + \int_{\Gamma_{\infty}}\nat(\ut^{-1}(\deb\ut\nabla\ut^{-1}+ \nabla\ut\deb\ut^{-1})\dd\ut ) \non\\
&=& -\nat\dd\int_{\Gamma_{\infty}}\ut^{-1}\nabla\deb \ut  + \int_{\Gamma_{\infty}}\nat\ut^{-1}\nabla\ut\nabla\ut^{-1}\dd\ut \ . \non
\eeq
Hence modulo coboundaries $\nat\dd(\cdot )$, the cyclic homology class of the anomaly is represented by
\beq
\Delta(\om,A) &=& \Phi\nat(\omt) - \frac{1}{2\pi i} \int_{\Gamma_{\infty}} \nat \nabla\At\, \omt \non\\
&\equiv& \Phi\nat(\ut^{-1}\dd\ut) - \frac{1}{4\pi i}\int_{\Gamma_{\infty}}\nat\ut^{-1}\nabla\ut\nabla\ut^{-1}\dd\ut \mod \nat\dd \non
\eeq
and coincides with $\sqrt{2\pi i}\, \ch(\rho_!(u))$ by Proposition \ref{pano}. Hence $\ch(\rho_!(u))$ appears as a cap-product of the Chern character $\ch_{\textup{dR}}(\ut)$ with the trace $\Phi$ and the Todd class. \\
Now let $e\in M_{\infty}(\Ac_0)$ be an idempotent. To calculate the pushforward of $[e]\in \Kt_0(\Ic\hotimes\Ac)$ we use Bott periodicity \cite{P5}. For notational simplicity suppose $e\in \Ac_0$. Let $S\Ac:=\Ac\hotimes\cinf(0,1)$ be the smooth suspension of $\Ac$. Under the Bott isomorphism $\Kt_0(\Ic\hotimes\Ac)\cong \Kt_1(\Ic\hotimes S\Ac)$ the class $[e]$ is represented by the invertible
$$
u= 1 + e\otimes (\beta -1) \in (\Ac_0\otimes\cinf(0,1))^+\ ,
$$
where $\beta  \in 1+ \cinf(0,1)$ is an invertible function with winding number $1$ (the Bott generator). The direct image $\rho_!(e)\in \Kt_0(\Ic\hotimes\Bc)$ then corresponds to $\rho_!(u)\in \Kt_1(\Ic\hotimes S\Bc)$ under Bott isomorphism. From \cite{P5} \S 4 the Chern character $\ch(\rho_!(e)) \in HP_0(\Bc)$ is related to $\ch(\rho_!(u)) \in HP_1(S\Bc)$ by the formula
$$
\ch(\rho_!(e)) \equiv \frac{1}{\sqrt{2\pi i}} \int_0^1 \ch(\rho_!(u)) \ \in \Th\Bc_{\nat} \ ,
$$
where $\int_0^1 : \Om^1\Th (S\Bc)_{\nat} \to \Th\Bc_{\nat}$ is the natural morphism $\Om^1\Th(\Bc\hotimes\cinf(0,1))\to \Om^1(\Th\Bc\hotimes\cinf(0,1))\to \Th\Bc \hotimes \Om^1(0,1)$ followed by integration of one-forms over the interval $[0,1]$. By construction one has $\ut-1 \in \Ac_0 \otimes \Th(\Bc\otimes\cinf(0,1))$. Let us still denote abusively by $\ut-1$ its image in $\Ac_0\otimes \Th\Bc \otimes \cinf(0,1)$. Let $s: \cinf(0,1)\to \Om^1(0,1)$ be the ordinary differential. Then we can write
$$
\ch(\rho_!(e)) \equiv \frac{1}{2\pi i} \int_0^1 \Big( \Phi\nat(\ut^{-1}s\ut) - \frac{1}{4\pi i}\int_{\Gamma_{\infty}}\nat\ut^{-1}\nabla\ut\nabla\ut^{-1}s\ut \Big) \ \in \Th\Bc_{\nat}
$$
Since it appears as a cap-product of $\ch_{\textup{dR}}(\ut)$ with the cyclic cocycles $\Phi$ and $\Td$, the class of the r.h.s. is homotopy invariant with respect to $\ut$. Then remark that $\ut$ is homotopic to the invertible
$$
v= 1 + \et \otimes (\beta -1) \ \in (\Ac_0\otimes\Th\Bc\otimes \cinf(0,1))^+\ ,
$$
because both $\ut$ and $v$ project to the same invertible in $(\Ac_0\otimes\Bc\otimes \cinf(0,1))^+$. Consequently the Chern character of $\rho_!(e)$ is represented by 
$$
\ch(\rho_!(e)) \equiv \frac{1}{2\pi i} \int_0^1 \Big( \Phi\nat(v^{-1}sv) - \frac{1}{4\pi i}\int_{\Gamma_{\infty}}\nat v^{-1}\nabla v \nabla v^{-1}sv \Big)\ .
$$
Let us calculate in the space $(\Omc^*(\Si)\cp G \otimes\Th\Bc)\otimes\Om^1(0,1)$:
$$
 v^{-1}sv = (1+\et(\beta^{-1}-1)) \et s\beta\ =\ \et \otimes \beta^{-1}s\beta \ , 
$$
\beq
 v^{-1}\nabla v \nabla v^{-1}sv &=& (1+\et(\beta^{-1}-1))\, \nabla\et (\beta-1)\, \nabla\et (\beta^{-1}-1)\, \et s\beta \non\\
&=& \et\nabla\et\nabla\et \otimes \beta^{-1}(\beta-1)(\beta^{-1}-1)s\beta \ . \non
\eeq
Then we obtain the desired formula for the Chern character
$$
\ch(\rho_!(e)) \equiv \Phi \nat (\et) - \frac{1}{2\pi i} \int_{\Gamma_{\infty}} \nat\, \et\nabla\et\nabla\et 
$$
from the integrals $\int_0^1 \beta^{-1}s\beta = 2\pi i$ and $\int_0^1 \beta^{-1}(\beta-1)(\beta^{-1}-1)s\beta = 4\pi i$.   \cqfd\\

\begin{remark}\textup{If $G$ acts without fixed points one has  $\Gamma_f = \varnothing$ and $\Gamma_{\infty} = \Gamma_0\subset \Si$. Hence the trace $\Phi(\Gamma)$ vanishes and the index theorem reduces to a cap-product with the Todd class 
\be
\Td(\Gamma) (a_0,a_1,a_2)=\int_{\Si} a_0\nabla a_1 \nabla a_2\ .
\ee
In \cite{P0} we already obtained this formula using the Hopf algebra of diffeomorphisms introduced by Connes and Moscovici \cite{CM95}. In the latter situation $G$ is a pseudogroup of conformal transformations whose action can be lifted to the bundle of K\"ahler metrics over $\Si$ by a Thom isomorphism, and the horizontal Dolbeault operator is combined with the vertical signature operator. Using characteristic classes for Hopf algebras as in \cite{CM98} we precisely obtained the above formula, up to an overall factor of 2 accounting for the contribution of the vertical signature operator. Note that the modular differential $\delta=[\d,D]$ is one of the generators of the Connes-Moscovici Hopf algebra. }
\end{remark}

\end{document}